\documentclass[12pt,leqno]{article}
\usepackage[official]{eurosym}
\usepackage{fullpage,amsmath,theorem}
\usepackage{epic,eepic}
\usepackage{amssymb,amsmath}
\usepackage{cases}

\usepackage{pictex}

\usepackage{color}
\usepackage{mathtools}
\usepackage{tikz-cd}
\usepackage{mathabx}

\usepackage{pdfpages}

\def\llvdash{{\|\hskip-2pt \raise 3pt\hbox{\vrule
height 0.25pt width 0.4cm}}}

\def\bet{{\beta}}

\def\calP{{\cal P}}

\def\l{\langle}
\def\r{\rangle}

\def\l{{\langle}}
\def\r{{\rangle}}

\def\calP{\mathcal P}

\def\oa{{\overline A^{\,\lower 7pt_{\hbox{$\scriptstyle\bet}}
\hbox{$\scriptstyle 0\tau$}}}}

\def\bet{\beta}

\def\llvdash{{\|\hskip-2pt \raise 3pt\hbox{\vrule height 0.25pt
width 0.4cm}}}

\newtheorem{theorem}{Theorem}[section]
\newtheorem{lemma}[theorem]{Lemma}

\newtheorem{proposition}[theorem]{Proposition}
{\theorembodyfont{\rmfamily}

{\theorembodyfont{\rmfamily}
\newtheorem{definition}[theorem]{Definition}}
{\theorembodyfont{\rmfamily}
\newtheorem{remark}[theorem]{Remark}}
{\theorembodyfont{\rmfamily}
\newtheorem{claim}{Claim}}
{\theorembodyfont{\rmfamily}
}
{\theorembodyfont{\rmfamily}
}
 \DeclareMathOperator{\dom}{dom}

\newcommand{\pr}{\medskip\noindent\textit{Proof}. }

\newcommand{\lusim}[1]{\smash{\underset{\raisebox{1.2pt}[0cm][0cm]{$\sim$}}
{{#1}}}}

\def\dom{{\rm dom}}

\def\llvdash{{\|\hskip-2pt \raise 3pt\hbox{\vrule
height 0.25pt width 0.15cm}}}

\def\Vdashbks{\hbox{$\Vdash\!\!\!\!{\raise2pt\hbox
{$\scriptscriptstyle\backslash$}}$}}

\newcommand{\R}{\mathbb{R}}
\newcommand{\PP}{\mathbb{P}}
\newcommand{\Q}{\mathbb{Q}}

\overfullrule=0pt
\begin{document}

\title{ The Galvin property at $\kappa^{++}$ and not at $\kappa^+$ }

\baselineskip=18pt
\author{  Moti Gitik  and Shachar Herpe
\footnote{ The work of the authors was partially supported  by ISF grant No. 882/22.}
 }

\date{\today}
\maketitle

\begin{abstract}
We construct a $\kappa-$complete ultrafilter $W$ over $\kappa$ such that $\neg$Gal$(\kappa, W, \kappa^+)$ and  Gal$(\kappa, W, \kappa^{++})$.
This answers a question of T. Benhamou and G. Goldberg \cite{B-Gol}.
\end{abstract}

\section{Introduction}
We deal here with the Galvin property of $\kappa$-complete ultrafilters over $\kappa$.

\begin{definition}[The Galvin property at $\lambda$] \label{galvin}
    Let $U$ be a filter over $\kappa$.
    Let $\lambda$ be a cardinal s.t.\ $\kappa < \lambda \leq 2^\kappa$.
    We say that $U$ satisfies the Galvin property at $\lambda$ (and denote it by Gal($\kappa,U,\lambda)$) iff for every $\{ A_\alpha | \alpha < \lambda \} \subseteq U$ there exists a sequence $\{ \alpha_i|i<\kappa\}\subseteq\lambda$ s.t.
    $$\underset{i<\kappa}{\bigcap} A_{\alpha_i}\in U.$$
\end{definition}

It was shown by F. Galvin  \cite{Bau-H-M} that if $2^{<\kappa}=\kappa$ then every normal filter $U$ over $\kappa$ satisfies Gal($\kappa,U,\kappa^+)$.
U. Abraham and S. Shelah \cite{A-She} showed consistency of $\neg$Gal$(\aleph_1, Cub_{\aleph_1}, \aleph_2)$.
The consistency of a negation of the Galvin property for $\kappa-$complete ultrafilters over $\kappa$ was first proved by T. Benhamou, S. Garti and S. Shelah in \cite{B-Gar-She}.
A supercompact cardinal was used for this. In \cite{B-G}, a different method was suggested and the initial assumption was reduced to a measurable.

Note that if Gal($\kappa,U,\lambda)$ holds for some ultrafilter $U$ and cardinal $\lambda$, then Gal($\kappa,U,\lambda')$ holds for every cardinal $\lambda'\geq \lambda$.
Equivalently, if $\lnot\text{Gal}(\kappa,U,\lambda)$ holds for some ultrafilter $U$ and cardinal $\lambda$, then $\lnot\text{Gal}(\kappa,U,\lambda')$ holds for every $\lambda'\leq\lambda$.
\\In all previously known examples,  it always was the case that either Gal($\kappa,U,\lambda)$ holds for every $\kappa < \lambda \leq 2^\kappa$ or $\lnot$ Gal($\kappa,U,\lambda)$ holds for every $\kappa < \lambda \leq 2^\kappa$.
This leads  to the natural question,  asked explicitly by T. Benhamou and G. Goldberg  \cite{B-Gol}:\\
\emph{Is  it consistent to have a $\kappa-$complete  ultrafilter $U$ over $\kappa$ such that} \ Gal($\kappa,U,2^\kappa$) and $\lnot$Gal($\kappa,U,\kappa^+)$?

The purpose of this paper is to provide an affirmative answer to this question.
Namely, we prove the following:

\begin{theorem}
    Assume GCH and that there is an elemantary embedding $j:V\rightarrow M$ with crit($j$)=$\kappa$ s.t. $(\kappa^{++})^M=\kappa^{++}$ and $\prescript{\kappa}{}{M}\subseteq M$.
    Then there is a cofinality preserving forcing extension $V^*$ such that $V^* \models 2^\kappa = \kappa^{++}$, and in $V^*$ there is a $\kappa$-complete ultrafilter $W$ over $\kappa$ s.t. $W$ does not satisfy Gal($\kappa,W,\kappa^+)$, but satisfies Gal($\kappa,W,\kappa^{++})$.
\end{theorem}

Notice that these assumptions are optimal, since in the extension  $\kappa$ is a measurable with $2^\kappa=\kappa^{++}$.

The proof follows the lines of Woodin  for blowing up the power of a measurable cardinal (as presented by Cummings in \cite{Cum})(it is possible to use \cite{Y} instead) and a construction from \cite{B-G}.
\\
The main part of the proof will be  to show that the constructed ultrafilter $W$, satisfies the Galvin property with respect to $\kappa^{++}$.
At the end of the paper some generalizations will be discussed.

We assume that the reader is familiar with the Woodin method for blowing up the power of a measurable cardinal using extenders, see Cummings \cite{Cum}.
Let us only state some basics.

\newpage
\section{Proof of the theorem}

Let $j:V\rightarrow M$ be as in the statement of the theorem.
We can derive a $(\kappa,\kappa^{++})$-extender $E$ from $j$.
Just set $E_a=\{ X\subseteq [\kappa]^{|a|}\mid a \in j(X)\},$ for every $a\in [\kappa^{++}]^{<\omega}$.
Further, let us identify $M$ with the ultrapower $M_E$ of $E$ and $j$ withits ultrapower embedding $j_E$.

\subsection{Blowing up the power of $\kappa$}
For each regular cardinal $\delta$, let Cohen($\delta,\delta^{++})$ be the Cohen forcing for adding $\delta^{++}-$Cohen functions to $\delta$. It consists of  partial functions from $\delta\times \delta^{++}$ to $\delta$ of cardinality less than $\delta$.
Let
$$\PP=\PP_{\kappa+1}=\langle \PP_\alpha,\lusim{\Q}_\beta|\alpha\leq \kappa+1,\beta\leq \kappa\rangle$$
be an Easton support iteration where
$$\Vdash_{\PP_\alpha} \Q_\alpha=\text{Cohen}(\check{\alpha},\check{\alpha}^{++})$$
when $\alpha$ is inaccessible and $\lusim{\Q}_\alpha$ is trivial otherwise.
\\Let $G*g_\kappa$ be generic for $\PP$, where $G$ is generic for $\PP_\kappa$ over $V$ and $g_\kappa$ is generic for $\Q_\kappa=Cohen(\kappa,\kappa^{++})$ over $V[G]$.
Denote  by $\langle f_{\kappa,\alpha} | \alpha<\kappa^{++} \rangle$  the generic Cohen functions added by   $g_\kappa$.

It is standard   to check that $V[G*g_\kappa] \models 2^\kappa=\kappa^{++}$.  Let us deal with measurability and extensions of elementary embeddings.
We will have to find a generic $H\in V[G*g_\kappa]$ for $j_E(\PP)$ s.t. $j_E''(G*g_\kappa)\subseteq H$ which will allow us to lift the embedding (as in Proposition 9.1 in \cite{Cum}).
In order to do that, we will use the projection of $E$ to its
 normal ultrafilter.

Let $U$ be the normal $\kappa$-complete ultrafilter on $\kappa$ derived from $j_E$:
$$U:=\{X\subseteq \kappa | \kappa\in j_E(X)\}.$$
Let $$ k:M_U\rightarrow M_E,\ k([f]_U)=j_E(f)(\kappa).$$
Denote $\lambda:=(\kappa^{++})^{M_U}$.

For convenience of the readers, let us state and proof below some well known claims.

\begin{claim}
\label{properties of k}
    The embedding $k$ is well defined, elementary, $j_E = k\circ j_U$ and crit($k)=\lambda$.
\end{claim}

\pr
$k$ is well defined since:
\begin{align*}
   & [f]_U=[g]_U\\
    \iff & \{\nu<\kappa|f(\nu)=g(\nu)\}\in U\\
    \iff & \kappa \in j_E(\{\nu<\kappa|f(\nu)=g(\nu)\})\\
    \iff & j_E(f)(\kappa)=j_E(g)(\kappa)\\
    \iff &k([f]_U)=k([g]_U).
    \end{align*}
$k$ is elementary since:
\begin{align*}
    &M_U \models \varphi([f_1]_U),\ldots,[f_n]_U) \\
    \iff & \{\alpha<\kappa|V\models \varphi(f_1(\alpha),\ldots,f_n(\alpha))\}\in U\\
    \iff & \kappa \in j_E(\{\alpha<\kappa|V\models \varphi(f_1(\alpha),\ldots,f_n(\alpha))\})\\
    \iff & M_E \models \varphi(j_E(f_1)(\kappa),\ldots,j_E(f_n)(\kappa))\\
    \iff &M_E \models \varphi(k([f]_U),\ldots,k([f_n]_U)).
\end{align*}
Now we have that $k(j_U(x))=k([C_x]_U)=j_E(C_x)(\kappa)=j_E(x)$.
Since $j_U(\kappa)$ is inaccessible in $M_U$, $\lambda < j_U(\kappa)$.
Since GCH holds, and since $j_U$ is an ultrafilter embedding, we get that $j_U(\kappa)<(2^\kappa)^+=\kappa^{++}$.
On the other hand, $k(\lambda)=(\kappa^{++})^{M_E}=\kappa^{++}>\lambda$, so crit($k)\leq\lambda$.
We claim that crit($k$) must be a cardinal (in $M_U$), since otherwise, if crit($k)=\beta > \alpha$ and $f:\alpha \leftrightarrow \beta$, then $k(f):\alpha\leftrightarrow k(\beta)$ and since
$$\forall \gamma<\alpha,\ k(f)(\gamma)=k(f)(k(\gamma))=k(f(\gamma))=f(\gamma),$$
Im($k(f))=\beta$ which is a contradiction.
Now we just need to show that $\kappa,\kappa^+\in \text{Im}(k)$.
Indeed, $k(\kappa)=k([Id]_U)=j_E(Id)(\kappa)=\kappa$.
Define
$$f:\kappa\rightarrow V \text{ by } f(\alpha)=\alpha^+.$$
Then $k([f]_U)=j_E(f)(\kappa)=\kappa^+$, which means that $\kappa^+\in \text{Im}(k)$, and thus crit($k)=\lambda$.
\\
$\square$

We will also use the following property of $k$:

\begin{claim}\label{width of k}
    The embedding $k$ has width $\leq \lambda$, i.e.,  every $x\in M_E$  is of the form $k(F)(a)$, for some $F\in M_U, a\in M_E$, where $M_U\models |\dom(F)|\leq \lambda$.

    \end{claim}

\pr
Let $x\in M_E$.
Then
$$x=j_E(f)(a)=k(j_U(f))(a)  \text{ where } f:[\kappa]^{|a|}\rightarrow V,\ a\in [\kappa^{++}]^{<\omega}.$$
Now, since $a\subseteq \kappa^{++}=k(\lambda)$, $a\in \text{dom}(k(j_U(f)\restriction \lambda))$.
Then
$$M_U\models |\text{dom}(j_U(f)\restriction \lambda)|=\lambda \text{ and }
 x=j_E(f)(a)=k(j_U(f))(a)=k(j_U(f)\restriction \lambda)(a).$$
\\
$\square$

Now from elementarity and $\kappa$-closure of $M_U$ and $M_E$ we get that the iterations $\PP$, $j_U(\PP)$ and $j_E(\PP)$ agree up to stage $\kappa$.

\begin{claim}
    $\PP=\PP_{\kappa+1}=(j_E(\PP))_{\kappa+1}$.
\end{claim}
\pr
    First of all, $(\kappa^{++})^{M_E}=\kappa^{++}$.
    Since $(j_E(\PP))_\kappa=\PP_\kappa$, we need to show that
    if $f$ is a partial function from $\kappa\times\kappa^{++}$, $f\in V[H]$ where $H$ is $\PP_\kappa-$generic and $|f|<\kappa$, then $f\in M_E[H]$.
    This is true since $\PP_\kappa$ is $\kappa^+$-c.c. and $\prescript{\kappa}{}M_E\subseteq M_E$, and then 
    $\prescript{\kappa}{}M_E[H]\subseteq M_E[H]$.
 \\
 $\square$

 \begin{claim}
     $(j_U(\PP))_{\kappa+1}=\PP_\kappa* \lusim{\Q}^*$ where $\lusim{\Q}^*=($Cohen($\kappa,\lambda))^{V[G]}$.
 \end{claim}

\pr
    This is because $\lambda=(\kappa^{++})^{M_U}$ and ${}^{\kappa}{}M_U[H]\subseteq M_U[H]$ for any $H$ which is $\PP_\kappa$ generic.
 \\
 $\square$

Let $g_0=g_\kappa\cap \Q_\kappa^*\subseteq M_U[G]$.
So $g_0$ is $\Q_\kappa^*$ generic over $V[G]$ (since every maximal antichain of $\Q_\kappa^*$ is also a maximal antichain of $\Q_\kappa$), and thus also over $M_U[G]$.

Let us first lift the embedding $j_E$ to $V[G]$.
In order to do so we will lift first $j_U$ to $V[G]$ and then we will lift $k$.
Denote $j_U(\PP)_{(\kappa,j_U(\kappa))}$ by $\R_0$ and $j_E(\PP)_{(\kappa,j_E(\kappa))}$ by $\R$.
So we have that
$$j_U(\PP_\kappa)=\PP_\kappa*\Q_\kappa^**\R_0 \text{ and } j_E(\PP_\kappa)=\PP_\kappa*\Q_\kappa*\R.$$
Since $\Q_\kappa^*$ is $\kappa^+-c.c.$ we get 
that
$$V[G*g_0] \models \prescript{\kappa}{}{M_U[G*g_0]}\subseteq M_U[G*g_0]$$
and thus $\R_0$ is $\kappa^+$-closed over $V[G*g_0]$.
Since
$$M_U[G*g_0] \models |\R_0|=j_U(\kappa) \land 2^{j_U(\kappa)}=j_U(\kappa)^+,$$
there are only $j_U(\kappa)^+$ antichains of $\R_0$ in $M_U[G*g_0]$.
But now in $V[G*g_0]$, $j_U(\kappa)^+$ has cardinality $\kappa^+$ 
and we can use Proposition 8.1 from \cite{Cum} to get a generic $H_0\in V[G*g_0]$ for $\R_0$ over $M_U[G*g_0]$.
In addition, since $k$ has width $\leq\lambda$ and $\R_0$ is $\lambda^+$-closed (in $M_U[G*g_0])$, we can use Proposition 15.1 in \cite{Cum} to get that $k''H_0$ generates a generic $H\in M_U[G*g_0*H_0]$ for $\R$ over $M_E[G*g_\kappa]$.

So we have the following commutative diagram:

\begin{center}
\begin{tikzcd}
    &V[G]\arrow[rr,"j_E"]\arrow[ddrr,"j_U"] & &M_E[G*g_\kappa * H]\\
    &&&\\
    & & &M_U[G*g_0*H_0]\arrow[uu,"k"]
\end{tikzcd}
\end{center}

Where we abuse the notation and denote the extension of $j_E$ to $V[G]$ still by $j_E$ (and similarly for $j_U$ and $k$).
It remains to find a generic
$$g_{j_E(\kappa)} \text{ for } j_E(\Q_\kappa) \text{ in } V[G*g_\kappa] \text{ s.t. }j_E'' g_\kappa \subseteq g_{j_E(\kappa)}.$$
To do so, let us first consider
$$\mathbb{S}_0:= j_U(\Q_\kappa)=(\text{Cohen}(j_U(\kappa),j_U(\kappa)^{++})^{M_U[G*g_0*H_0]}.$$
We will force with $\mathbb{S}_0$ over $V[G*g_\kappa]$ and find the generic there.
So we need to show that doing so will not change cofinalities.

\begin{claim}
    $\mathbb{S}_0$ is $\kappa^+$-closed and $\kappa^{++}$-Knaster in $V [G * g_0]$.
\end{claim}

\pr
    $\mathbb{S}_0$ is $\kappa^+$-closed because $M_U[G*g_0*H_0] \models \mathbb{S}_0 \text{ is } j_U(\kappa)$ closed.
    So if $\langle p_\alpha : \alpha<\kappa \rangle$ is an
    increasing sequence of conditions in $V[G*g_0]$,
    $$V [G * g_0] \models \prescript{\kappa}{} M_U[G * g_0] \subseteq M_U[G * g_0],$$
    and thus
    $$\langle p_\alpha : \alpha<\kappa \rangle \in M_U[G * g_0].$$
    Now since $\kappa<j_U(\kappa)$, it has an upper bound.

    To see that $\mathbb{S}_0$ is $\kappa^{++}$-Knaster, let
    $$\langle p_\alpha : \alpha < \kappa^{++}\rangle\in V[G*g_0]$$
    be a sequence of conditions from $\mathbb{S}_0$, and let
    $$p_\alpha = j_U(f_\alpha)(\kappa) \text{ where }f_\alpha : \kappa \rightarrow \Q_\kappa, f_\alpha \in V[G].$$
    Define
    $$A_\alpha=\langle \{\beta\}\times\text{dom}(f_\alpha(\beta))| \beta < \kappa \rangle,\ \forall \alpha<\kappa^{++}.$$
    Then $|A_\alpha|=\kappa$ and $V[G] \models (\kappa^+)^\kappa=\kappa^+$, so we can use the $\Delta$-system theorem to get
    \begin{align*}
        &B\subseteq \kappa^{++}, |B|=\kappa^{++}\text{ and }\langle r_i | i<\kappa \rangle \text{ s.t. }\\
        &\forall \alpha_1,\alpha_2 \in B, \forall \beta<\kappa, \text{dom}(f_{\alpha_1}(\beta))\cap\text{dom}(f_{\alpha_2}(\beta))=r_\beta.
    \end{align*}
    Now since $\forall i<\kappa,\ |r_i|=\delta_i<\kappa$, there are only $$\underset{i<\kappa}{\Pi} 2^{\delta_i}=\kappa$$ options for functions $f$ with dom($f(\beta))=r_\beta$ and thus we can find $B'\subseteq B$ of size
    $\kappa^{++}$ s.t. $\forall \alpha_1, \alpha_2 \in B'$, $f_{\alpha_1}$ and
    $f_{\alpha_2}$ are pointwise compatible.
    Then $\langle p_\alpha| \alpha\in B' \rangle$ are pairwise compatible.
 \\
 $\square$

\begin{claim}
    $\mathbb{S}_0$ is $(\kappa^+,\infty)$-distributive and $\kappa^{++}$-c.c. in $V[G * g_\kappa]$.
\end{claim}

\pr
    Notice that we can use a poset $\Q$ which is isomorphic to $\Q_\kappa$ to extend $V[G*g_0]$ to $V[G*g]$.
    Now the poset $\Q_\kappa$ is $\kappa^+$-c.c. in $V[G * g_0]$ and 
     since $\Q$ must also be $\kappa^+$-c.c., $\mathbb{S}_0$ is $(\kappa^+,\infty)$-distributive in $V [G* g_\kappa]$.
     $\mathbb{S}_0 \times \Q$ is $\kappa^{++}$-c.c. in $V [G * g_0]$ and so 
      we get that $\mathbb{S}_0$ is $\kappa^{++}$-c.c. in $V[G*g_\kappa]$.
 \\
 $\square$

So we get that $\mathbb{S}_0$ preserves cofinalities, and we can now force with it over $V[G*g_\kappa]$.
Denote the generic by $h_0$.
It is also clear that
$$V[G*g_\kappa*h_0] \models 2^\kappa =\kappa^{++}.$$
Now we want to transfer $h_0$ along $k$ to get a generic $g$ for $\mathbb{S}:=k(\mathbb{S}_0)=j_E(\Q_\kappa)$.
This is again possible by Proposition 15.1 in \cite{Cum}, since $k$ has width $\leq\lambda$ and
$$M_U[G*g_0*H_0] \models \mathbb{S}_0 \text{ is } \lambda^+\text{-closed}.$$
Now $g$ does not have to satisfy that $j_E''g_\kappa \subseteq g$.
So let us change its conditions to satisfy it, in a way that does not ruin the genericity.
In order to build the ultrafilter $W$ later, we will also need a few additional changes.
So for every $p\in j_E(\Q_\kappa)$, let $p^*$ be defined as follows:
$\text{dom}(p)=\text{dom}(p^*)$ and for every $\langle \gamma,\alpha \rangle \in \text{dom}(p)$,
\begin{align*}
    p^*(\langle \gamma,\alpha \rangle) =
    \begin{cases}
        f_{\kappa,\beta}(\gamma) & \gamma<\kappa \land \alpha = j_E(\beta)\\
        2\cdot \beta & \gamma=\kappa \land \alpha=j_E(\beta)<j_E(\kappa^+)\\
        \beta & \gamma=\kappa \land \alpha=j_E(\kappa^+ + \beta)\\
        p(\langle \gamma, \alpha \rangle) & \text{else}
    \end{cases}
\end{align*}

\begin{claim}
$p^*\in M_E[G*g_\kappa*H]$.
\end{claim}

\pr
    Since $\prescript{\kappa}{}{M_E[G*g_\kappa*H]}\subseteq M_E[G*g_\kappa*H]$, it is enough to show that $p$ was only changed in $\kappa$ many places.
    Let $a \in [\kappa^{++}]^{<\omega}$ such that $j_E(f)(a) = p$, where $f : \kappa^{|a|} \rightarrow \Q_\kappa$.
    By elementarity, if
    $$\langle \alpha, j_E(\beta) \rangle \in \kappa \times j_E''\kappa^{++} \cap (\text{dom}(p)),$$
    there is
    $$x \in \kappa^{|a|} \text{ such that }\langle \alpha, \beta \rangle \in \text{dom}(f(x)).$$
It follows that $|(\kappa \times j_E''\kappa^{++})\cap \text{dom}(p)| \leq \kappa$.
Moreover, $|\{\kappa\}\times j_E''\kappa^{++}\cap\text{dom}(p)| \leq \kappa$, since otherwise there would be some $\alpha < \kappa^{++}$ such that
$$\text{cf}(\alpha) = \kappa^+ \text{ and sup}\{j_E(\beta) | \langle \kappa, j_E(\beta)\rangle \in \text{dom}(p)\} = j_E(\alpha).$$
But $|\text{dom}(p)|^{M_E[G*g_\kappa*H]} < j_E(\kappa)$ and $(\text{cf}(j_E(\alpha)))^{M_E[G*g_\kappa*H]} = j_E(\kappa^+)$ which is a contradiction.
 \\
 $\square$

Now define $g_{j_E(\kappa)}:= \{p^*|p\in g\}$.

\begin{lemma}
\label{genericity after change}
The filter $g_{j_E(\kappa)}$ is Cohen$(j_E(\kappa), j_E(\kappa)^{++})^{M_E[G*g_\kappa*H]}$-generic over $M_E[G * g_\kappa * H]$.
\end{lemma}
\pr
    Let us work in $M_E[G*g_\kappa*H]$.
    Let $D$ be open dense.
    Define $D^*$ to consist of all conditions $p\in \text{Cohen}(j_E(\kappa), j_E(\kappa)^{++})$ such that
    \begin{align*}
        \forall q. \text{ if }(\text{dom}(q) = \text{dom}(p) \land |\{x | p(x) \neq q(x)\}| \leq \kappa\\
        \land (p(x)\neq q(x)\rightarrow q(x)<\kappa^{++})) \text{ then } q \in D.
    \end{align*}
    Let us show that $D^*$ is dense.
    Let  $p\in \text{Cohen}(j_E(\kappa), j_E(\kappa)^{++})$ and enumerate by $\langle q_r | r < \theta\rangle$ all the conditions $q$ such that
    $$\text{dom}(q) = \text{dom}(p) \land |\{x | p(x) \neq q(x)\}| \leq\kappa \land (p(x)\neq q(x)\rightarrow q(x)<\kappa^{++}).$$
    Let $|\text{dom}(p)|=\alpha<j_E(\kappa)$.
    Then $\theta\leq [\alpha]^{\kappa}\cdot(\kappa^{++})^\kappa\leq 2^\alpha\cdot 2^{\kappa^{++}}< j_E(\kappa)$ since $j_E(\kappa)$ is inaccessible.
    We define inductively an increasing sequence $\langle p_r | r < \theta\rangle$.
    Define $p_0 = p$, and suppose that $p_r$ is defined.
    Let $p'_{r+1} := q_r \cup p_r \restriction (\text{dom}(p_r) \setminus \text{dom}(p))$.
    Choose $p'_{r+1}\leq t_{r+1}\in D$ (which exists by density) and set $p_{r+1} = p_r \restriction \text{dom}(p) \cup t_{r+1} \restriction (\text{dom}(t_{r+1}) \setminus \text{dom}(p))$.
    Then $p_r \leq p_{r+1}$.
    When $r$ is limit, take just the union (it is possible since $\mathbb{Q}_{j_E(\kappa)}$ is $j_E(\kappa)$-closed).
    Let
    $$p':= \underset{r<\theta}{\bigcup} p_r.$$
    Then $p'$ has the property that any $\kappa$ many changes of $p'$ from the domain of $p$ into $\kappa^{++}$ belongs to $D$.
    Namely for any $q\in \text{Cohen}(j_E(\kappa), j_E(\kappa)^{++})$, if $\text{dom}(q) = \text{dom}(p')$,
    $$q \restriction (\text{dom}(p') \setminus \text{dom}(p)) = p'\restriction (\text{dom}(p') \setminus \text{dom}(p)),$$
    $$|\{x \in \text{dom}(p) | p'(x) \neq q(x)\}| \leq\kappa$$ and
    $$p'(x)\neq q(x)\rightarrow q(x)<\kappa^{++},$$
    then $q\in D$.
    This is because $q \restriction \text{dom}(p) = q_r$ for some r, therefore $q \geq t_{r+1} \in D$ and $D$ was open.
    Now we define inductively $\langle p^{(r)} | r < \kappa^{+}\rangle,\ p^{(0)} = p$, at limit steps we take union, and at successor steps we take
    $$p^{(r+1)} := (p^{(r)})'.$$
    We claim that
    $$p\leq p_* := \underset{r<\kappa^{+}}{\bigcup} p^{(r)} \in D^*.$$
    First note that $|p_*| < j_E(\kappa)$ so $p_*\in \text{Cohen}(j_E(\kappa),j_E(\kappa)^{++})$.
    Let $q$ be any condition with dom$(q) = \text{dom}(p_*)$ and let
    $$I = \{x \in \text{dom}(p_*) | q(x) \neq p_*(x)\},$$
    and suppose that
    $$|I| \leq \kappa \land \forall x\in I. q(x)<\kappa^{++}.$$
    Since
    $$\text{dom}(p_*) = \underset{r<\kappa^+}{\bigcup} \text{dom}(p^{(r)})$$
    and dom$(p^{(r)})$ is $\subseteq$-increasing, there is $j < \kappa^+$ such that $I \subseteq \text{dom}(p^{(j)})$.
    The condition $q\restriction \text{dom}(p^{(j)})$ is enumerated in the construction of $p^{(j+1)}$, hence $q \restriction \text{dom}(p^{(j+1)}) \in D$, and since $D$ is open, $q \in D$.
    This means that $p_* \in D^*$.
    Finally, by genericity of $g$, we can find $p \in D^* \cap g$.
    By definition, $p^*\in g_{j_E(\kappa)}$, and since
    $$\text{dom}(p^*) = \text{dom}(p), |\{x | p(x) \neq p^*(x)\}| \leq \kappa \text{ and }p^*(x)\neq p(x) \rightarrow p^*(x)<\kappa^{++},$$
    it follows that $p^*\in  D$.
 \\
 $\square$

Now we can lift $j_E$ to $j_E^*:V[G*g_\kappa]\rightarrow M_E[G*g_\kappa*H*g_{j_E(\kappa)}]$.
But since we obtained $g_{j_E(\kappa)}$ in $V[G*g_\kappa*h_0]$, $j_E^*$ is not necessarily definable in $V[G*g_\kappa]$.
We need also to find a generic for $j_E(\mathbb{S}_0)$.
This is again possible from Proposition~15.1 in \cite{Cum}, since $j_E$ has width $\leq \kappa$ and $V[G*g_\kappa] \models \mathbb{S}_0$ is $(\kappa^+,\infty)$-distributive.
Denote the generic by $h_1$ to get
$$j_E^* : V [G * g_\kappa * h_0] \rightarrow M_E[G * g_\kappa * H * g_{j_E(\kappa)} * h_1].$$
Again, with some abuse of notation, we denote the lift of $j_E^*$ to $V[G*g_\kappa*h_0]$ again by $j_E^*$.

\subsection{Constructing $W$ and proving $\lnot$Gal$(\kappa,W, \kappa^+)$}

Let $\kappa_1:= j_E(\kappa)$ and denote by $\langle f_{\kappa_1,\alpha} | \alpha<\kappa_1^{++} \rangle$ the Cohen functions of $g_{\kappa_1}$.
Denote:
\begin{align*}
G_{\kappa_1}:= G*g_\kappa * H , \  G_{\kappa_1+1}:= G_{\kappa_1}*g_{\kappa_1}, \ G_{\kappa+1}:=G*g_\kappa.
\end{align*}
In addition,  denote $V[G_{\kappa+1}*h_0]$ by $V^*$ and $M_E[G_{\kappa_1+1}*h_1]$ by $M_E^*$.
Let $U^*$ be the normal ultrafilter derived from $j_E^*$:
$$U^* = \{X \subseteq \kappa | \kappa \in j_E^*(X)\}.$$

\begin{lemma}
    $j_E^*=j_{U^*}$.
\end{lemma}

\pr
Let $k:M_{U^*}\rightarrow M_E^*$ be the elementary embedding defined by $k([f]_{U^*})=j_E^*(f)(\kappa)$.
Let us show that $k$ is onto.
So let $x=(\lusim{x})_{G_{\kappa_1+1}*h_1}\in M_E^*$.
Since $\lusim{x}\in M_E$, there are some
$$g \in V , a = \{\alpha_1,\ldots, \alpha_r\} \in [\kappa^{++}]^{<\omega} \text{ such that }j_E(g)(a) = \lusim{x}.$$
Define in $V^*$ the function
$$g^*(\alpha) = (g(\{f_{\kappa,\kappa^++\alpha_1}(\alpha), ..., f_{\kappa,\kappa^++\alpha_r}(\alpha)\}))_{G_{\kappa+1}*h_0}.$$
Then,
\begin{align*}
    k([g^*]_{U^*}) &= j_E^*(g^*)(\kappa)\\
    &= (j_E^*(g)(\{j_E^*(f_{\kappa,\kappa^++\alpha_1})(\kappa),...,j_E^*(f_{\kappa,\kappa^++\alpha_r})(\kappa)\}))_{G_{\kappa_1+1}*h_1}\\
    &= (j_E(g)(\{f_{\kappa_1,j_E(\kappa^++\alpha_1)}(\kappa),...,f_{\kappa_1,j_E(\kappa^++\alpha_r)}(\kappa)\}))_{G_{\kappa_1+1}*h_1}\\
    &= (j_E(g)(\{\alpha_1,...,\alpha_r\}))_{G_{\kappa_1+1}*h_1}\\
    &=(j_E(g)(a))_{G_{\kappa_1+1}*h_1}\\
    &=(\lusim{x})_{G_{\kappa_1+1}*h_1}=x.
\end{align*}

Notice that the fourth equality follows from the construction of $g_{\kappa_1}$, since we have:
$$\forall \alpha<\kappa^{++},\ f_{\kappa_1,j_E(\kappa^+ +\alpha)}(\kappa)=\alpha.$$
 \\
 $\square$

Now let us consider the second ultrapower (of $V$) by $E$, i.e. Ult($M_E,j_{j_E(E)})$.
To facilitate notations, let us denote:
$$M_1 := M_E,\ M_2 := M_{j_E(E)}^{M_1},\ j_1:= j_E,\ i:=j_{j_E(E)}, \ j_2:= i\circ j_1.$$
Let also $\kappa_2 := j_2(\kappa)=i(\kappa_1)$.
So we have the following diagram:

\begin{center}
\begin{tikzcd}
    &V\arrow[rr,"j_2"]\arrow[ddrr,"j_1"] & &M_2(=M_{j_E(E)})\\
    &&&\\
    & & &M_1(=M_E)\arrow[uu,"i"]
\end{tikzcd}
\end{center}

Now we wish to extend $i$ to $M_1^*$.
By taking the second ultrapower by $U^*$, by elementarity, we get that $M_{(U^*)^2}\simeq M_{j_{U^*}(U^*)}$ is a generic extension of $M_2$.
This is because
$$V^* \models (j_{U^*}\restriction_V=j_E) \text{ and thus } M_1^*\models (j_{j_{U^*}(U^*)}\restriction_{M_1} = j_{j_E(E)}).$$
In addition, since $M_1^*$ is a generic extension of $M_1$ by $j_1(\PP*\mathbb{S}_0)$, $M_{(U^*)^2}$ is a generic extension of $M_2$ by $j_2(\PP*\mathbb{S}_0)$.
So we can extend $i$ to $M_1^*$.
Denote the generic that we get for $j_2(\PP*\mathbb{S}_0)$ by $G_{\kappa_2}*g_{\kappa_2}'*h_2'$.
So now we have the following diagram:

\begin{center}
\begin{tikzcd}
    &V[G*g_\kappa*h_0]\arrow[rr,"j_{(U^*)^2}"]\arrow[ddrr,"j_{U^*}"] & &M_2[G_{\kappa_2}*g'_{\kappa_2}*h'_2]\\
    &&&\\
    & & &M_1[G_{\kappa_1}*g_{\kappa_1}*h_1]\arrow[uu,"j_{j_{U^*}(U^*)}"]
\end{tikzcd}
\end{center}

In order to construct $W$, we need to change $g'_{\kappa_2}$.
So for every $p\in$ Cohen($\kappa_2,\kappa_2^{++})^{M_2[G_{\kappa_2}]}$, define $p^*$
 with dom($p^*)$=dom$(p)$ and for every $\langle \gamma,\alpha \rangle \in$ dom($p$) set
\begin{align*}
    p^*(\langle \gamma, \alpha \rangle) = \begin{cases}
        2 \cdot \beta + 1 & \gamma=\kappa_1 \land \alpha=i(\beta), \alpha \in j_2''\kappa^+\\
        p(\langle \gamma, \alpha \rangle) & else
    \end{cases}
\end{align*}

Note that from the definition of $g_{\kappa_1}$, we will also get that $p^*(\langle\kappa,i(\beta)\rangle)=2 \cdot \beta$ when  $i(\beta)\notin j_2''\kappa^+$ and $\beta<\kappa_1^+$ and $p^*(\langle\kappa,i(\kappa_1^++\beta)\rangle)=\beta$.
We also have that $|\pi_2(\text{dom}(p))\cap j_2''\kappa^+|\leq \kappa$ ) (where $\pi_2$ is the projection to the second coordinate).
This follows since $j_2''\kappa^+$ is unbounded in $\kappa_2^+$, and $\pi_2(\text{dom}(p))$ is bounded in $\kappa_2^+$.
So we get that $p^*\in M_1^*$.\\
To make sure that it is also in $M_2[G_{\kappa_2}]$, we need to check that we only changed $p$ in $\kappa_1$ many places.
In other words, we want to show that
$$M_1^* \models |\pi_2(\text{dom}(p))\cap i''\kappa_1^+|\leq \kappa_1.$$
But this is clear, since $i''\kappa_1^+$ is unbounded in $\kappa_2^+$ whereas $|p|<\kappa_2$.

Hence $p^* \in M_2[G_{\kappa_2}]$ since we have only changed $p$ at $\kappa_1$ many values and
$$\prescript{\kappa_1}{}M_2[G_{\kappa_2}]\cap M_1[G_{\kappa_1}] \subseteq M_2[G_{\kappa_2}].$$
Setting $g_{\kappa_2}=\{p^*|p\in g_{\kappa_2}'\}$, we get that it is still generic for Cohen$(\kappa_2,\kappa_2^{++})^{M_2[G_{\kappa_2}]}$ over $M_2[G_{\kappa_2}]$.
This can be shown in the same way as in Lemma \ref{genericity after change}.

Since $i''g_{\kappa_1}\subseteq g_{\kappa_2}'$, and we did not change $g_{\kappa_2}'$ below $\kappa_1$, we can lift $i$ to $i^*$:
$$i^*:M_1^*\rightarrow M_2[G_{\kappa_2}*g_{\kappa_2}*h_2]=:M_2^*.$$

Then the embedding $j_2 : V \rightarrow M_2$ extends to
$$i^*\circ j_1^* = j_2^*:  V^* \rightarrow M_2^*.$$

Denote the functions of $g_{\kappa_2}$ by $\{f_{\kappa_2,\alpha}|\alpha<\kappa_2^{++}\}$.
So we get in particular that
\begin{enumerate}
    \item $f_{\kappa_2,i(\alpha)}(\kappa_1)$ is odd if $\alpha\in j_1''\kappa^+$.
    \item $f_{\kappa_2,i(\alpha)}(\kappa_1)$ is even if $\alpha \notin j_1''\kappa^+$.
\end{enumerate}

Define now $$W = \{X \subseteq \kappa^2 |   (\kappa,\kappa_1) \in j_2^*(X)\}.$$

We will show that $W$ is the desired ultrafilter.

\begin{remark}
    Notice that $W$ is Rudin-Kiesler equivalent to some ultrafilter over $\kappa$, but it will be more convinient to work with $W$.
\end{remark}
To see this, let us first observe some important information about $W$:

\begin{lemma}\label{W not kappa+ glavin}
$j_W=j_2^*$ and $ [id]_W=(\kappa,\kappa_1)$.
\end{lemma}

\pr
    Set $k:M_W\rightarrow M_2^*$ by $k([f]_W)=j_2^*(f)(\kappa,\kappa_1)$.
    To see that $k$ is onto, let $(\lusim{x})_{G_{\kappa_2+1}*h_2}\in M_2^*$.
    Since $i$ is an ultrapower by a $(\kappa_1,\kappa_1^{++})$-extender, there are
    $$f\in M_1 \text{ and } a=\{\alpha_1,...,\alpha_r\}\in [\kappa_1^{++}]^{<\omega}, \ \text{ s.t. } \lusim{x}=i(f)(a).$$
    Define
    $$M_1^* \ni g:\kappa_1 \rightarrow \text{On},\ g(\gamma)=(f(f_{\kappa_1,\kappa_1^++\alpha_1}(\gamma),...,f_{\kappa_1,\kappa_1^++\alpha_r}(\gamma))_{G_{\kappa_1+1}*h_1}.$$
    Then, by the construction of $g_{\kappa_2}$,
    \begin{align*}
        i^*(g)(\kappa_1)&=(i(f)(f_{\kappa_2,i(\kappa_1^++\alpha_1)}(\kappa_1),...,f_{\kappa_2,i(\kappa_1^++\alpha_r)}(\kappa_1)))_{G_{\kappa_2+1}*h_2}\\
        &=(i(f)(\alpha_1,...,\alpha_r))_{G_{\kappa_2+1}*h_2}\\
        &=(i(f)(a))_{G_{\kappa_2+1}*h_2}\\
        &=(\lusim{x})_{G_{\kappa_2+1}*h_2}=x.
    \end{align*}
    Now since $M_1^*$ is an ultrapower by a normal ultrafilter, we get that there is some
    $$q:\kappa\rightarrow V^* \text{ s.t. } g=j_1^*(q)(\kappa)$$
    and then
    $$x=i^*(j_1^*(q)(\kappa))(\kappa_1)=j_2^*(q)(\kappa)(\kappa_1).$$
    In particular, we can assume that $q(\beta): \kappa \rightarrow V^*$ for every $\beta<\kappa$.
    Now define $f^*:\kappa^2\rightarrow V^*$ by
    $$f^*(\alpha,\beta)=q(\alpha)(\beta).$$
    So we get that
    \begin{align*}
        k([f^*]_W)&=j_2^*(f^*)(\kappa,\kappa_1)\\
        &=j_2^*(q)(\kappa)(\kappa_1)\\
        &=x.
    \end{align*}
    In particular we get that
    $$[id]_W=k([id]_W)=j_2^*(id)(\kappa,\kappa_1)=(\kappa,\kappa_1).$$
 \\
 $\square$

Finally, for every $\alpha<\kappa^+$, let us define
$$ A_\alpha:=\{(\gamma,\beta)\in\kappa^2|f_{\kappa,\alpha}(\beta) \text{ is odd}\}.$$
We changed $g'_{\kappa_2}$ in such a way that these sets will be in $W$ and will witness that $W$ does not satisfy the Galvin property at $\kappa^+$:
\begin{lemma}
    $\{ A_\alpha|\alpha<\kappa^+\}\subseteq W$ and this sequence witnesses  $\lnot\text{Gal}(\kappa,W,\kappa^+)$.
\end{lemma}
\pr Repeat the corresponding argument from \cite{B-G}.

\subsection{Proving Gal$(\kappa,W,\kappa^{++})$}
\subsubsection{Some preparations}

We arrive now to \textbf{the main issue of this thesis,} which is to show that $\text{Gal}(\kappa,W,\kappa^{++})$ holds.
Let $\{ B_\alpha | \alpha<\kappa^{++} \} \subseteq W$.

For each $\alpha<\kappa^{++}$, we  choose a nice name for $B_\alpha$ (in V[G]),
$$\lusim{B}_\alpha'=\{\{i\}\times A^{\alpha}_i|i<\kappa\}.$$
Since Cohen($\kappa,\kappa^{++})^{V[G]}$ is $\kappa^+$-c.c., there are only $\kappa$ conditions in $A_i^\alpha$.
Each condition has a support of size $<\kappa$. Set
$$a_\alpha = {\bigcup}_{i<\kappa} {\bigcup}_{p\in A^{\alpha}_i} \pi_2''\text{dom}(p)\subseteq \kappa^{++}.\footnote{ $\pi_2$ denotes here the projection of pairs to the second coordinate.}$$
Then $|a_\alpha|\leq\kappa$.
Now we  use the $\Delta$-system lemma to get a subset $A$ of $\kappa^{++}$ of size $\kappa^{++}$ and
$$a\subseteq \kappa^{++} \text{ s.t. }\forall \alpha\neq\beta\in A, \ a_\alpha\cap a_\beta = a.$$
Since
$$\forall \alpha,\beta\in A, \text{ if } \alpha\neq\beta \text{ then }(a_\alpha\cap a_\beta \cap \kappa^+) \setminus a=\emptyset,$$
there are only $\kappa^+$ non-empty options for $(a_\alpha \setminus a)\cap \kappa^+$.
Pick any $\langle B_{\alpha_i} | i < \kappa^{++} \rangle$ which satisfies that
$$a_{\alpha_i}\cap a_{\alpha_j}=a \text{ and }(a_{\alpha_i} \setminus a)\cap \kappa^+=\emptyset, \ \forall i\neq j<\kappa^{++}.$$
Without loss of generality we can assume that $\alpha_i=i,\ \forall i<\kappa^{++}$.
Since we would like to work in $V$ (and not in $V[G]$), let $\lusim{B}_\alpha$ be a name in $V$ s.t. $(\lusim{B}_\alpha)_G=\lusim{B}_\alpha'$ and $(\lusim{B}_\alpha)_{G*g}=B_\alpha$.

For each $\alpha<\kappa^{++}$ choose $q^\alpha\in j_2(\PP)_{<\kappa_2}$, $p^\alpha,r^\alpha\in \text{Cohen}(\kappa_2,\kappa_2^{++})^{M_2[G_{\kappa_2}]}$ and names $\lusim{p}^\alpha,\ \lusim{r}^\alpha$ (in $V[G_{\kappa_2}]$) s.t.
\begin{enumerate}
    \item $(\lusim{r}^\alpha)_{G_{\kappa_2}}=r^\alpha,\ (\lusim{p}^\alpha)_{G_{\kappa_2}}=p^\alpha$
    \item $q^\alpha\in G_{\kappa_2}$
    \item $r^\alpha\in g_{\kappa_2}'\cap g_{\kappa_2}$
    \item $\text{dom}(r^\alpha)\subseteq \kappa_2\times \kappa_2^{++}\setminus \kappa_2^+  $
    \item $p^\alpha\in g_{\kappa_2}$
    \item $\text{dom}(p^\alpha)\subseteq \kappa_2\times\kappa_2^+$
    \item ${q^\alpha}^\frown (\lusim{p}^\alpha\cup \lusim{r}^\alpha) \Vdash (\check{\kappa},\check{\kappa_1})\in j_2(\lusim{B}_\alpha)$
\end{enumerate}

Now we want to stabilize the $p^\alpha$'s.
Let us look at
$$Z:=\text{Cohen}(j_U(\kappa),j_U(\kappa^+))^{M_U[G*g_0*H_0]}.$$
Recall that we had the following diagram:
\begin{center}
\begin{tikzcd}
    &V[G]\arrow[rr,"j_E"]\arrow[ddrr,"j_U"] & &M_E[G*g_\kappa * H]\\
    &&&\\
    & & &M_U[G*g_0*H_0]\arrow[uu,"k"]
\end{tikzcd}
\end{center}

and that $h_0$ was a generic for $j_U(\Q_\kappa)$ over $M_U[G*g_0*H_0]$.

So
$$M_U[G*g_0*H_0]\models Z \text{ has } j_U(\kappa^+) \text{ open dense sets}$$
So in $V^*$ we can enumerate them by $\langle D_\xi|\xi<\kappa^+\rangle$.
Now we can build an increasing sequence of conditions
$$\langle r_\xi|\xi<\kappa^+\rangle\subseteq Z\text{, s.t. }r_\xi\in h_0\cap D_\xi.$$
Since we had that $k''h_0$ generates $g_{\kappa_1}'$ we get that
$$\langle k(r_\xi)|\xi<\kappa^+\rangle \text{ generates }g_{\kappa_1}'\restriction \text{Cohen}(\kappa_1,\kappa_1^+)^{M_E[G*g_\kappa*H]}=g_{\kappa_1}'\restriction k(Z).$$
Now look at
$$j_1(\langle D_\xi|\xi<\kappa^+\rangle)=\langle D^1_\xi|\xi<\kappa_1^+\rangle,\ j_1(\langle r_\xi|\xi<\kappa^+\rangle)=\langle r^1_\xi|\xi<\kappa_1^+\rangle.$$
Recall that $M_2=M_{j_E(E)}^{M_E}$, so by elementarity, $\langle D^1_\xi|\xi<\kappa_1^+\rangle$ is an enumeration of the open dense sets of $j_1(Z)$ and $\langle r^1_\xi|\xi<\kappa_1^+\rangle$ is an increasing sequence of conditions in $j_1(Z)$ s.t.
$$r^1_\xi\in D^1_\xi\cap j_1^*(h_0).$$
Then we get again that
$$\langle j_E(k)(r_\xi^1)|\xi<\kappa_1^+\rangle \text{ generates }g_{\kappa_2}\restriction \text{Cohen}(\kappa_2,\kappa_2^+)^{M_2[G_{\kappa_2}]}.$$

But since the sequence is increasing, and $j_E''\kappa^+$ is unbounded in $\kappa_1^+$, we can actually have that $$\langle j_E(k)(r_{j_E(\xi)}^1)|\xi<\kappa^+\rangle \text{ generates }g_{\kappa_2}\restriction \text{Cohen}(\kappa_2,\kappa_2^+)^{M_2[G_{\kappa_2}]}.$$
This means that for every $\alpha<\kappa^{++}$ there is some
$$\xi_\alpha<\kappa^{+}\text{ s.t. }p^\alpha\subseteq j_E(k)(r_{j_E(\xi_\alpha)}^1).$$
But then we have that for $\kappa^{++}$ many $\alpha$'s it is the same $\xi^*$, and we can just use $j_E(k)(r_{j_E(\xi^*)}^1)$ instead of all these $p^\alpha$'s.
Let us simply denote it by $p$.
Notice that technically the above sequences do not generate $g_{\kappa_1}\restriction \kappa_1^+$ and $g_{\kappa_2}\restriction \kappa_2^+$, but the generics we had before the change.
But this is not a problem since we just need to stabilize the domain (because they all come from the same generic).

Now we will actually need to work with $t^\alpha,m^\alpha,s^\alpha$ where
$${t^\alpha}^{\frown} {\lusim{m}^\alpha}^{\frown} \lusim{s}^\alpha={q^\alpha}^\frown \lusim{r}^\alpha \text{ and }t^\alpha\in j_2(\PP)_{\kappa_1}\text{, }m_\alpha\in \Q_{\kappa_1}\text{ and }s^\alpha\in j_2(\PP)_{(\kappa_1,\kappa_2]}.$$

Let $h_{t^\alpha},h_{m^\alpha},h_{s^\alpha},h_p^\alpha$ represent them.
We want them all to be in $V$, so they depend on some generator $\rho_\alpha<\kappa^{++}$, i.e.,
\begin{align*}
    j_2(h_{t^\alpha})(\rho_\alpha,j_1(\rho_\alpha))&=t^\alpha,\\
    j_2(h_{m^\alpha})(\rho_\alpha,j_1(\rho_\alpha))&=\lusim{m}^\alpha,\\
    j_2(h_{s^\alpha})(\rho_\alpha,j_1(\rho_\alpha))&=\lusim{s}^\alpha.\\
    j_2(h_p^\alpha)(\rho_\alpha,j_1(\rho_\alpha))&=\lusim{p}
\end{align*}

In addition, let $h_p\in V^*$ represent $\lusim{p}$, so we have that $j_2(h_p)(\kappa,\kappa_1)=\lusim{p}$.
Let  $\pi_\alpha$ denotes a projection of $E_{\rho_\alpha}$ to $U$, where
$$E_{\rho_\alpha}=\{X\subseteq \kappa| \rho_\alpha\in j_1(X)\}.$$
So $j_1(\pi_\alpha)(\rho_\alpha)=\kappa$ and it follows that $j_2(\pi_\alpha)(j_1(\rho_\alpha))=\kappa_1$.
Notice that the set
$$\{(\nu_1,\nu_2)\in\kappa^2|h_p(\pi_\alpha(\nu_1),\pi_\alpha(\nu_2))=h_p^\alpha(\nu_1,\nu_2)\}\in (E_{\rho_\alpha}^*)^2,$$
so we can assume that this is true for every $(\nu_1,\nu_2)\in \kappa^2$ (otherwise change $h_p^\alpha$ on a small set).
Now since
$$j_2(h_{t^\alpha})(\rho_\alpha,j_1(\rho_\alpha))\in \PP_{\kappa_1}=\PP_{j_2(\pi_\alpha)(j_1(\rho_\alpha))}$$
we can assume wlog that
$$\forall (\nu_1,\nu_2)\in \kappa^2.h_{t^\alpha}(\nu_1,\nu_2)\in \PP_{\pi_\alpha(\nu_2)}.$$
Similarly we can assume that $h_{m^\alpha}(\nu_1,\nu_2)\in \lusim{Q}_{\pi_\alpha(\nu_2)}$ and $h_{s^\alpha}(\nu_1,\nu_2)\in \lusim{\PP}_{(\pi_\alpha(\nu_2),\kappa]}$ for all $(\nu_1,\nu_2)\in \kappa^2$.

Another general fact we'll need is that on a set in $U^*$, $\pi_\alpha\circ f_{\kappa,\kappa^++\rho_\alpha}=id$.
This is because we changed the generics s.t. $f_{\kappa_1,j_1(\kappa^++\rho_\alpha)}(\kappa)=\rho_\alpha$ and thus
$$j_1^*(\pi_\alpha\circ f_{\kappa,\kappa^++\rho_\alpha})(\kappa)=j_1(\pi_\alpha)(j_1^*(f_{\kappa,\kappa^++\rho_\alpha})(\kappa))=j_1^*(\pi_\alpha)(\rho_\alpha)=\kappa.$$

In particular, since $U^*$ is normal, we get that $\pi_\alpha$ is $1-1$ on a set in $U^*$.
Denote this set by $D_\alpha$ and notice that $D_\alpha\in U^*\cap V[G*R]$ where
$$    R:=\{f_{\kappa,\alpha}|\kappa^+\leq \alpha\}\subseteq g_\kappa.$$

Let $\alpha<\kappa^{++}$. Define the following sets:

\begin{align*}
&A_\alpha:= \\
&\{\vec{\nu}\in \kappa^2|h_{t^\alpha}(\vec{\nu})^\frown h_{m^\alpha}(\vec{\nu})^\frown (h_{s^\alpha}(\vec{\nu})\cup h_p^\alpha(\vec{\nu})) \Vdash (\check{\pi_\alpha(\nu_1)},\check{\pi_\alpha(\nu_2)})\in \lusim{B}_\alpha\}
\end{align*}

and
\begin{align*}
    A_{t^\alpha}&:=\{(\nu_1,\nu_2)\in \kappa^2|h_{t^\alpha}(\nu_1,\nu_2)\in G_{\pi_\alpha(\nu_2)}\},\\
    A_{m^\alpha}&:=\{(\nu_1,\nu_2)\in \kappa^2|(h_{\lusim{m}^\alpha}(\nu_1,\nu_2))_{G_{\pi_\alpha(\nu_2)}}\in g_{\pi_\alpha(\nu_2)}\}, \\
    A_{s^\alpha}&:=\{(\nu_1,\nu_2)\in\kappa^2|(h_{\lusim{s}^\alpha}(\nu_1,\nu_2))_{G_{\pi_\alpha(\nu_2)}*g_{\pi_\alpha(\nu_2)}}\in G_{(\pi_\alpha(\nu_2),\kappa)}*R\},\\
    A_p&:=\{(\nu_1,\nu_2)\in\kappa^2|(h_p(\nu_1,\nu_2))_G\in g_\kappa\}.
\end{align*}
So we have that $A_{t^\alpha},A_{m^\alpha},A_{s^\alpha}\in V[G*R]$ and
\begin{align*}
    B_\alpha&\supseteq \pi_\alpha''(A_\alpha\cap A_{t^\alpha}\cap A_{m^\alpha}\cap A_{s^\alpha})\cap A_p\\
    &=\{(\pi_\alpha(\nu_1),\pi_\alpha(\nu_2))\in A_p|(\nu_1,\nu_2)\in A_\alpha\cap A_{t^\alpha}\cap A_{m^\alpha}\cap A_{s^\alpha}\}.
\end{align*}

We also have that $A_p\in W$ since $(j_2^*(h_p)(\kappa,\kappa_1))_{G_{\kappa_2}}=p\in g_{\kappa_2}$.
In addition $A_\alpha\in (E_{\rho_\alpha})^2$ and similarly, the sets $A_{t^\alpha},A_{m^\alpha},A_{s^\alpha}$ are in $(E_{\rho_\alpha}^*)^2$ where
$$E_{\rho_\alpha}^*=\{X\subseteq \kappa| \rho_\alpha\in j_1^*(X)\}.$$

We come now to the main point:
\begin{lemma}\label{the idea}
    Suppose  that   there is a sequence
    $$\l X_{A_\alpha},X_{A_{t^\alpha}},X_{A_{m^\alpha}},X_{A_{s^\alpha}}\mid \alpha<\kappa^{++}\r$$ in  $ V[G*R]$ which consists of sets in $U^*$
    such that
    $$[X_{A_\alpha}]^2\subseteq \pi_\alpha''A_\alpha,\ [X_{A_{t^\alpha}}]^2\subseteq \pi_\alpha''A_{t^\alpha},\ [X_{A_{m^\alpha}}]^2\subseteq \pi_\alpha''A_{m^\alpha},\ [X_{A_{s^\alpha}}]^2\subseteq \pi_\alpha''A_{s^\alpha},$$ for every $\alpha<\kappa^{++}$.\\
    Then there is some $I\subseteq \kappa^{++}$ s.t.
    $$|I|=\kappa\text{ and }{\bigcap}_{\alpha\in I} B_\alpha\in W.$$
\end{lemma}
\pr
    Let
    $$X_\alpha:=X_{A_\alpha}\cap X_{A_{t^\alpha}}\cap X_{A_{m^\alpha}}\cap X_{A_{s^\alpha}}\cap D_\alpha\in U^*\cap V[G*R].$$
    \begin{claim}
    \label{the sets are there}
    For every $X\in U^*\cap V[G*R]$, $[X]^2\in W$.
\end{claim}
\pr
    We have that
    $$V^*=V[G*g_\kappa*h_0]=V[G*R*(g_\kappa\setminus R)*h_0].$$
    Denote the embedding $j_2^*\restriction_{V[G*R]}$ by $j_2^1$.
    So we get that
    $$j_2^1:V[G*R]\rightarrow M_2[G_{\kappa_2}*R_2] \text{ where } R_2=j_2^*(R)=\{f_{\kappa_2,\alpha}|\alpha\geq \kappa_2^+\}.$$
    Now let us denote $j_{(U^*)^2}$ by $j_2'$.
    So if we show that $j_2'\restriction _{V[G*R]}=j_2^1$, we will get that
    $$\forall X\in V[G*R],\ j_2'(X)=j_2^*(X).$$

    But this is clear since $$j'_2(R)=\{f'_{\kappa_2,\alpha}| \alpha\geq \kappa_2^+\}=\{f_{\kappa_2,\alpha}| \alpha\geq \kappa_2^+\}=R_2.$$
    Now if $X\in U^*\cap V[G*R]$ we get that
\begin{align*}
\kappa\in j_1^*(X)&\implies \kappa\in j_2'(X)\land j_2^*(X)\land \kappa_1\in j_2'(X)=j_2^*(X)\\
&\implies (\kappa,\kappa_1)\in j_2^*([X]^2)\\
&\implies [X]^2\in W.
\end{align*}
 \\
 $\square$

    We will repeat now the original proof of Galvin in order to get a set
    $$I\subseteq \kappa^+ \text{ such that } |I|=\kappa \text{ and }\underset{\alpha\in I}{\bigcap}X_\alpha\in U^*.$$
    The proof is done inside $V[G*R]$. In particular, $I$ and $\underset{\alpha\in I}{\bigcap}X_\alpha$ will be in $V[G*R]$. Then, by the claim above,    $[\underset{\alpha\in I}{\bigcap}X_\alpha]^2\in W$. So, we will obtain  the following:
    \begin{align*}
        {\bigcap}_{\alpha\in I}B_\alpha&\supseteq \underset{\alpha\in I}{\bigcap}(  \pi_\alpha''(A_\alpha\cap A_{t^\alpha}\cap A_{m^\alpha}\cap A_{s^\alpha}))\cap A_p\\
        &\supseteq {\bigcap}_{\alpha\in I}[X_\alpha]^2\cap A_p \\
        &=[{\bigcap}_{\alpha\in I}X_\alpha]^2\cap A_p\in W
    \end{align*}

    For the reader convenience,  let us go through  the Galvin proof to see that we can stay within $V[G*R]$.
    For every $\alpha<\kappa^+,\ \xi<\kappa$, let
    $$H_{\alpha,\xi}=\{\beta<\kappa^+|X_\alpha\cap \xi=X_\beta\cap \xi \}.$$
    Since we have that $\langle X_\alpha|\alpha<\kappa^+\rangle\in V[G*R]$, clearly $H_{\alpha,\xi}\in V[G*R]$.
    \begin{claim}
    There is some $\alpha^*<\kappa^+$ s.t. for every $\xi<\kappa$, $|H_{\alpha^*,\xi}|=\kappa^+$.
    \end{claim}
    \pr
    Assume otherwise.
    Then for every $\alpha<\kappa^+$ let $\xi_\alpha<\kappa$ s.t. $|H_{\alpha,\xi_\alpha}|<\kappa^+$.
    By the pigeonhole principle, there is some $\xi^*<\kappa$ and a set $A\subseteq \kappa^+$ with $|A|=\kappa^+$ s.t. $\forall\alpha\in A$, $\xi_\alpha=\xi^*$.
    But since $\kappa$ is measurable, $2^{\xi^*}<\kappa$ so there is some $A'\subseteq A$ of size $\kappa^+$  and $C\subseteq \xi^*$ s.t. $\forall \alpha\in A'$, $X_\alpha\cap \xi^*=C$.
    And this is a contradiction since if $\alpha\in A'$ then
    $A'\subseteq H_{\alpha,\xi^*} \text{ and } |A'|=\kappa^+>|H_{\alpha,\xi^*}|.$
     \\
 $\square$

    Let $\alpha^*<\kappa^+$ as in the claim.
    So we can define an increasing sequence $\langle \alpha_i|i<\kappa\rangle\in V[G*R]$ s.t. $\alpha_i\in H_{\alpha^*,i+1}$.
    \begin{claim}
        $X_{\alpha^*}\cap \underset{i<\kappa}{\triangle}X_{\alpha_i}\subseteq \underset{i<\kappa}{\bigcap}X_{\alpha_i}$.
    \end{claim}
    \pr
        Let
        $$\beta\in X_{\alpha^*}\cap \underset{i<\kappa}{\triangle}X_{\alpha_i} \text{ and } j<\kappa.$$
        If $j<\beta$, then $\beta\in \underset{i<\kappa}{\triangle}X_{\alpha_i}$ and thus $\beta\in X_{\alpha_j}$.
        If $ \beta\leq j$, then $\beta\in X_\alpha^*$ and $X_{\alpha^*}\cap (j+1)=X_{\alpha_j}\cap (j+1)$ which means that $\beta\in X_{\alpha_j}$.
     \\
 $\square$

    It is clear from the proof above that $\underset{i<\kappa}{\bigcap}X_{\alpha_i}\in V[G*R]$, and we get that $\underset{i<\kappa}{\bigcap}B_{\alpha_i}\in W$ as needed.
 \\
 $\square$

\subsubsection{Finding the one dimensional sets}
It is now left to find the sets $X_{A_\alpha},X_{A_{t^\alpha}},X_{A_{m^\alpha}},X_{A_{s^\alpha}}$.

We will need the following wellknown general claim:
\begin{claim} \label{one dimensional}
Let $\rho<\kappa^{++}$ and let $\pi$ represent $\kappa$ with respect to $E_\rho$.
Then for every $Y\in E_\rho^2$, there is some $X\in E_\rho$ s.t. $[X]^{2*}\subseteq Y$, where
$$[X]^{2*}:=\{(\alpha,\beta)\in X^2|\alpha<\pi(\beta)\}.$$
\end{claim}
\pr
    Set
    $$Y^*:=\{(\alpha,\beta)\in Y|\alpha<\pi(\beta)\},\ Y^*_1:=\{\alpha<\kappa|\exists \beta.(\alpha,\beta)\in Y^*\}.$$

    First we claim that $Y^*\in E_\rho^2$.
    This is because the set $[\kappa]^{2*}\in E_\rho^2$, since
    $$\forall\alpha<\kappa.\alpha<j_E(\pi)(\rho)=\kappa$$
    and thus
    $$\{\alpha<\kappa|\{\beta<\kappa|\alpha<\pi(\beta)\}\in E_\rho\}=\kappa\in E_\rho.$$
Now for every $\alpha\in Y_1^*$ set
$$Y_\alpha^*:=\{\beta<\kappa|(\alpha,\beta)\in Y^*\}.$$
So there is a set $Z\in E_\rho$ s.t. $\forall \alpha\in Z.Y_\alpha^*\in E_\rho$.
Set
$$X:=Z\cap{\triangle^*}_{\alpha\in Z}Y_\alpha^*$$
where
$${\triangle^*}_{\alpha\in Z}Y_\alpha^*=\{\nu<\kappa| \forall \alpha\in Z\cap \pi(\nu).\nu\in Y_\alpha^*\}.$$

Now $X$ is as desired since if $(\alpha,\beta)\in X^2$ and $\alpha<\pi(\beta)$, then $\beta\in Y_\alpha^*$ and thus $(\alpha,\beta)\in Y^*\subseteq Y$.
 \\
 $\square$

Let us start  with $A_\alpha$:
\begin{claim}
    There is a set $X_{A_\alpha}\in U^*\cap V[G*R]$ s.t. $[X_{A_\alpha}]^2\subseteq \pi_\alpha'' A_\alpha$.
\end{claim}
\pr
As in claim \ref{one dimensional} we can get an $X_\alpha\in E_{\rho_\alpha}$ s.t. $[X_\alpha]^{2*}\subseteq A_\alpha$.
Now we define
$$A'_\alpha:=\{\nu<\kappa | f_{\kappa,\kappa^++\rho_\alpha}(\nu)\in X_\alpha\}.$$
Now since
$$j_1^*(f_{\kappa,\kappa^++\rho_\alpha})(\kappa)=\rho_\alpha\in j_1(X_\alpha),$$
we have that $A_\alpha'\in U^*$.
In addition we define
$$C_\alpha:=\{\nu<\kappa|f_{\kappa,\kappa^++\rho_\alpha}(\nu)<\nu^{++} \text{ and }\nu \text{ is inaccessible}\}.$$

So $C_\alpha\in U^*$ since $\kappa$ is inaccessible and $\rho_\alpha<\kappa^{++}$.

Now we have that
$$ \pi_\alpha'' A_\alpha \supseteq \pi_\alpha''[X_\alpha]^{2*}$$
and we claim that

$$f_{\kappa,\kappa^++\rho_\alpha}''[(A_\alpha'\cap C_\alpha\cap D_\alpha)]^2\subseteq[X_\alpha]^{2*}.$$
To see this let $(\nu_1,\nu_2)\in [A'_\alpha\cap C_\alpha\cap D_\alpha]^2$.
Then
$$f_{\kappa,\kappa^++\rho_\alpha}(\nu_1)<\nu_1^{++}<\nu_2 \text{ and } \pi_\alpha(f_{\kappa,\kappa^++\rho_\alpha}(\nu_2))=\nu_2>f_{\kappa,\kappa^++\rho_\alpha}(\nu_1)$$
and thus
$$(f_{\kappa,\kappa^++\rho_\alpha}(\nu_1),f_{\kappa,\kappa^++\rho_\alpha}(\nu_2))\in[X_\alpha]^{2*}.$$
But in total we have that
$$\pi_\alpha''A_\alpha \supseteq \pi_\alpha''[X_\alpha]^{2*}\supseteq\pi_\alpha '' f_{\kappa,\kappa^++\rho_\alpha}'' [A_\alpha'\cap C_\alpha\cap D_\alpha]^2=[A_\alpha'\cap C_\alpha\cap D_\alpha]^2.$$
So $X_{A_\alpha}:=A_\alpha'\cap C_\alpha\cap D_\alpha\in U^*\cap V[G*R]$ is as desired.
 \\
 $\square$

Now, let us deal with $A_{t^\alpha}$:

\begin{claim}
    There is a set $X_{A_{t^\alpha}}\in U^*\cap V[G*R]$ s.t. $[X_{A_{t^\alpha}}]^2\subseteq \pi_\alpha'' A_{t^\alpha}$.
\end{claim}
\pr
Notice that $t^\alpha\in M_1$ so we can choose $h_{t^\alpha}$ s.t. $j_1(h_{t^\alpha})(\rho_\alpha)=t^\alpha$.
So we get that
$$t^\alpha=i''t^\alpha=i(t^\alpha)=i(j_1(h_{t^\alpha})(\rho_\alpha))=j_2(h_{t^\alpha})(\rho_\alpha)$$
which means that $h_{t^\alpha}$ does not depend on the second coordinate.
In this case we actually have that
$$A_{t^\alpha}=\{(\nu_1,\nu_2)\in \kappa^2|h_{t^\alpha}(\nu_1)\in G\}=\{\nu<\kappa|h_{t^\alpha}(\nu)\in G\}\times\kappa.$$

Let us define
$$A_{t^\alpha,1}:=\{\nu<\kappa|h_{t^\alpha}(\nu)\in G\},\ A_{t^\alpha,2}:=\{\nu<\kappa|f_{\kappa,\kappa^++\rho_\alpha}(\nu)\in A_{t^\alpha,1}\}.$$
So $A_{t^\alpha,1}\in E_{\rho_\alpha}^*$ and thus $A_{t^\alpha,2}\in U^*$.
So clearly
$$f_{\kappa,\kappa^++\rho_\alpha}''[A_{t^\alpha,2}\cap D_\alpha]^2\subseteq A_{t^\alpha}$$
and then,
$$\pi_\alpha''f_{\kappa,\kappa^++\rho_\alpha}''[A_{t^\alpha,2}\cap D_\alpha]^2=[A_{t^\alpha,2}\cap D_\alpha]^2\subseteq \pi_\alpha''A_{t^\alpha}$$
and $X_{A_{t^\alpha}}:=A_{t^\alpha,2}\cap D_\alpha\in U^*\cap V[G*R]$ is as desired.
 \\
 $\square$

Now, let us deal with $A_{s^\alpha}$:
\begin{claim}
    There is a set $X_{A_{s^\alpha}}\in U^*\cap V[G*R]$ s.t. $[X_{A_{s^\alpha}}]^2\subseteq \pi_\alpha'' A_{s^\alpha}$.
\end{claim}
We have $j_2(h_{s^\alpha})(\rho_\alpha,j_1(\rho_\alpha))=\lusim{s}^\alpha$.
Define
$$h^*_{s^\alpha}(\nu):\kappa\rightarrow \PP \text{ by } h_{s^\alpha}^*={\sup}_{\nu'<\pi_\alpha(\nu)}(h_{s^\alpha}(\nu',\nu)).$$
It is well defined since for all $\nu'<\pi_\alpha(\nu)$, $h_{s^\alpha}(\nu',\nu)\in \lusim{\PP}_{(\pi_\alpha(\nu),\kappa]}$ which is $\pi_\alpha(\nu)^+$-closed.
Then we have that
$$\lusim{s}^\alpha=j_2(h_{s^\alpha})(\rho_\alpha,j_1(\rho_\alpha))\leq j_2(h_{s^\alpha}^*)(j_1(\rho_\alpha)).$$
Set
\begin{align*}
C_{s^\alpha}&:=\{(\nu_1,\nu_2)\in [\kappa]^2|h_{s^\alpha}(\nu_1,\nu_2)\leq h_{s^\alpha}^*(\nu_2)\},\\
A_{s^\alpha}^*&:=\{\nu<\kappa|(h_{s^\alpha}^*(\nu))_{G_{\pi_\alpha(\nu)}*g_{\pi_\alpha(\nu)}}\in G_{(\pi_\alpha(\nu),\kappa)}*R\}.
\end{align*}
So we have that $C_{s^\alpha}\in E_{\rho_\alpha}^2$ and $A_{s^\alpha}^*\in E_{\rho_\alpha}^*$.
Let $X_{s^\alpha}\in E_{\rho_\alpha}$ s.t. $[X_{s^\alpha}]^{2*}\subseteq C_{s^\alpha}$ (as in claim \ref{one dimensional}).
Now denote $X_{s^\alpha}\cap A_{s^\alpha}^*$ by $F_{s^\alpha}$ and notice that $F_{s^\alpha}\in E_{\rho_\alpha}^*$.
Now we continue as before: let
\begin{align*}
A_{s^\alpha}'&:=\{\nu<\kappa|f_{\kappa,\kappa^++\rho_\alpha}(\nu)\in F_{s^\alpha}\},\\
C_\alpha^s&:=\{\nu<\kappa|f_{\kappa,\kappa^++\rho_\alpha}(\nu)<\nu^{++}\land \nu \text{ is inaccessible}\}.
\end{align*}
Clearly $A_{s^\alpha}',C_\alpha^s\in U^*\cap V[G*R]$.
Then $f_{\kappa,\kappa^++\rho_\alpha}''[C_\alpha^s\cap A_{s^\alpha}'\cap D_\alpha]^2\subseteq [F_{s^\alpha}]^{2*}$.
We also have that $[F_{s^\alpha}]^{2*}\subseteq (\kappa\times A_{s^\alpha}^*)\cap C_{s^\alpha}$, since if $\nu_1,\nu_2\in F_{s^\alpha}$ and $\nu_1<\pi_\alpha(\nu_2)$ then $(\nu_1,\nu_2)\in [X_{s^\alpha}]^{2*}\subseteq C_{s^\alpha}$ and clearly $\nu_2\in A_{s^\alpha}^*$.
Now it is only left to notice that $(\kappa\times A_{s^\alpha}^*)\cap C_{s^\alpha}\subseteq A_{s^\alpha}$ and apply $\pi_\alpha''$ to see that $X_{A_{s^\alpha}}:=C_\alpha^s\cap A_{s^\alpha}'\cap D_\alpha$ is as desired.
 \\
 $\square$

Now what is left is to deal with $A_{m^\alpha}$, which will also require a bit more work.
\begin{claim}
    There is a set $X_{A_{m^\alpha}}\in U^*\cap V[G*R]$ s.t. $[X_{A_{m^\alpha}}]^2\subseteq \pi_\alpha'' A_{m^\alpha}$.
\end{claim}
\pr
Notice that $m_\alpha\in \Q_{\kappa_1}=j_1^*(\Q_\kappa)$, so we can work in $M_1$.
We would like to stabilize $m^\alpha\restriction \text{Cohen}(\kappa_1,\kappa_1^+)$.
Recall that we defined
$$Z:=\text{Cohen}(j_U(\kappa),j_U(\kappa^+))^{M_U[G*g_0*H_0]}$$
and that $h_0$ was the generic for $\text{Cohen}(j_U(\kappa),j_U(\kappa^{++}))^{M_U[G*g_0*H_0]}$ (see p. 31 for more details).
Since $g_{\kappa_1}$ is generated by
$$\langle k(x)|x\in h_0\restriction Z\rangle \text{ and } V\models |h_0\restriction Z|=\kappa^+,$$
there is some condition $m\in \text{Cohen}(\kappa_1,\kappa_1^+)$ and some $A\subseteq\kappa^{++}$ with $|A|=\kappa^{++}$ s.t. for every $\alpha\in A$, $m^\alpha\restriction\text{Cohen}(\kappa_1,\kappa_1^+)\subseteq m$.
To avoid more indexes, we assume that $A=\kappa^{++}$ and replace $m^\alpha\restriction\text{Cohen}(\kappa_1,\kappa_1^+)$ by $m$.
As in the case of $p$, we can have some $h_m$ s.t. $j_2(h_m)(\kappa,\kappa_1)=\lusim{m}$.
So we have now that $m^\alpha=m\cup m_\alpha^*$ where $m_\alpha^*:=m^\alpha\restriction [\kappa_1^+,\kappa_1^{++})$.
Let $h_{m_\alpha^*}$ represent $\lusim{m}$ with respect to $E_{\rho_\alpha}$.
In particular we can define
\begin{align*}
    A_m&:= \{(\nu_1,\nu_2)\in [\kappa]^2|(h_m(\nu_1,\nu_2))_{G_{\nu_2}}\in g_{\nu_2}\}\\
    A_{m^\alpha}^*&:=\{(\nu_1,\nu_2)\in [\kappa]^2|(h_{m_\alpha^*}(\nu_1,\nu_2))_{G_{\pi_\alpha(\nu_2)}}\in g_{\pi_\alpha(\nu_2)}\restriction[\pi_\alpha(\nu_2)^+,\pi_\alpha(\nu_2)^{++})\}
\end{align*}
Notice that $A_m\in W$ and $A_m\cap \pi_\alpha''A_{m^\alpha}^*\subseteq \pi_\alpha''A_{m^\alpha}$.
So we only need to deal with $A_{m^\alpha}^*$ since $A_m$ does not depend on $\alpha$.

Notice that
$$(\kappa_1^{++})^{M_2}=(\kappa_1^{++})^{M_1}=j_1(\kappa^{++})=\bigcup j_1''\kappa^{++}$$
(and let us denote it simply by $\kappa_1^{++}$).
Consider the set
$$C_2:=\{\tau<\kappa_1^{++}|\exists h\in V[G*R].\tau=j_2^*(h)(j_1(\rho_\alpha))\}.$$
Then $C_2$ is unbounded in $\kappa_1^+$ and $\kappa_1^{++}$, since $\forall \beta<\kappa^{++}$,
\begin{align*}
j_2^*(f_{\kappa,\kappa^++\beta}\circ \pi_\alpha)(j_1(\rho_\alpha))&=f_{\kappa_2,\kappa_2^++j_2(\beta)}(j_2(\pi_\alpha)(j_1(\rho_\alpha))\\
&=f_{\kappa_2,\kappa_2^++j_2(\beta)}(\kappa_1)\\
&=j_1(\beta).
\end{align*}

Since $\sup(j_1''\kappa^+)=\kappa_1^+$ and $\sup(j_1''\kappa^{++})=\kappa_1^{++}$ we get in particular that $C_2\cap j_1''\kappa^+$ is unbounded in $\kappa_1^+$ and $C_2\cap j_1''\kappa^{++}$ is unbounded in $\kappa_1^{++}$.

Now fix some $\alpha<\kappa^{++}$ and let us deal with it.
Since $\pi_2''\text{dom}(m_\alpha^*)$ is bounded in $\kappa_1^{++}$, we can choose some $\eta_\alpha^2\in C_2\cap j_1''\kappa^{++}$ above its supremum.
Let $\sigma_\alpha^2: \kappa_1^+\leftrightarrow \eta_\alpha^2\setminus\kappa_1^+$ be the least $M_2$ bijection (which is also the least $M_1$ bijection).

Now $|\text{dom}(m_\alpha^*)|<\kappa_1$, so $(\sigma_\alpha^2)^{-1}(\pi_2''\text{dom}(m_\alpha^*))$ is bounded in $\kappa_1^+$.
Choose some $\eta_\alpha^1\in C_2\cap j_1''\kappa^+$ above it and $\sigma_\alpha^1:\kappa_1\leftrightarrow \eta_\alpha^1$ to be a least $M_2$ bijection (which is again the least $M_1$ bijection).

So
$$ \pi_2''\text{dom}(m_\alpha^*)\subseteq \sigma^2_\alpha{''}\sigma^1_\alpha{''}\kappa_1=\sigma^2_\alpha{''}\eta_\alpha^1.$$
Now since $|\text{dom}(m_\alpha^*)|<\kappa_1$, there is some $\eta_\alpha^0<\kappa_1$ s.t.
$$\pi_2''\text{dom}(m_\alpha^*)\subseteq\sigma_\alpha^2{''}\sigma_\alpha^1{''}\eta_\alpha^0.$$
Set
$$x_\alpha:=\sigma^2_\alpha{''}\sigma^1_\alpha{''}\kappa_1,\ x_\alpha':=\sigma_\alpha^2{''}\sigma_\alpha^1{''}\eta_\alpha^0.$$
Notice that since $\sigma_\alpha^1,\sigma_\alpha^2\in M_1$, so are $x_\alpha,x_\alpha'$ and we can look at $i(x_\alpha),\ i(x_\alpha')$ (where $i:M_1\rightarrow M_2$ is the embedding derived from $j_E(E)$).
Let $\varphi_\alpha$ be the bijection from $x_\alpha$ to $i''x_\alpha$.
So there are $\tilde{\eta}_\alpha^1,\tilde{\eta}_\alpha^2<\kappa^{++}$ and $h_{\eta_\alpha^1},h_{\eta_\alpha^2}\in V$ s.t.
$$\eta_\alpha^1=j_1(\tilde{\eta}_\alpha^1)=j_2(h_{\eta_\alpha^1})(j_1(\rho_\alpha)) \text{ and }\eta_\alpha^2=j_1(\tilde{\eta}_\alpha^2)=j_2(h_{\eta_\alpha^2})(j_1(\rho_\alpha))$$
We can also get that $\sigma_\alpha^1,\sigma_\alpha^2$ depend only on $h_{\eta_\alpha^1},h_{\eta_\alpha^2}$, since we can define $h_{\sigma_\alpha^2}(\nu)$ to be the minimal bijection from $(\pi_\alpha(\nu))^+$ to $h_{\eta_\alpha^2}(\nu)\setminus (\pi_\alpha(\nu))^+$ and $h_{\sigma_\alpha^1}(\nu)$ to be the minimal bijection from $\pi_\alpha(\nu)$ to $h_{\eta_\alpha^1}(\nu)$.
Then we can let
$$h_{x_\alpha}(\nu)=h_{\sigma_\alpha^2}(\nu){''}h_{\sigma_\alpha^1}(\nu){''}\pi_\alpha(\nu).$$
As for $\varphi_\alpha$, notice that $i(\sigma_\alpha^2)$ is the $M_2$ minimal bijection from $j_2(\tilde{\eta}_\alpha^2)\setminus \kappa_2^+$ to $\kappa_2^+$ and  $i(\sigma_\alpha^1)$ is the $M_2$ minimal bijection from $\kappa_2$ to $j_2(\tilde{\eta}_\alpha^1)$.
This is because $\sigma_\alpha^i$ is also $M_1$ minimal.
So we get that
$$\varphi_\alpha(\xi)=i(\xi)=i(\sigma_\alpha^2(\sigma_\alpha^1(\tau_\xi)))=i(\sigma_\alpha^2)(i(\sigma_\alpha^1)(\tau_\xi))$$ (for some $\tau_\xi<\kappa_1$).

So we can define $h_{i(\sigma_\alpha^2)}(\nu)$ to be the minimal bijection from $\kappa^+$ to $\tilde{\eta}_\alpha^2\setminus \kappa^+$ and $h_{i(\sigma_\alpha^1)}(\nu)$ to be the minimal bijection from $\kappa$ to $\tilde{\eta}_\alpha^1$ (and thus they do not depend on $\nu$).
Mark them as $\tilde{\sigma}_\alpha^2$ and $\tilde{\sigma}_\alpha^1$.
Let us define $h_{\varphi_\alpha}$ s.t.
$$h_{\varphi_\alpha}(\nu)(\xi)=\tilde{\sigma}_\alpha^2(\tilde{\sigma}_\alpha^1(\tau_\xi)),\ \forall \xi\in h_{x_\alpha}(\nu)$$
where $\xi=h_{\sigma_\alpha^2}(\nu)(h_{\sigma_\alpha^1}(\nu)(\tau_\xi))$.
Notice that the $\tau_\xi$'s goes over $\pi_\alpha(\nu)$.

Now define:
$$E^\alpha:= \{\nu<\kappa|\forall \xi\in h_{x_\alpha}(\nu),\ f_{\pi_\alpha(\nu),\xi}=f_{\kappa,h_{\varphi_\alpha(\nu)(\xi)}}\restriction\pi_\alpha(\nu)\}$$
and recall that we defined
$$E_{\rho_\alpha}^*=\{X\subseteq \kappa| \rho_\alpha\in j_1^*(X)\}.$$
Notice that $\kappa^+\leq h_{\varphi_\alpha(\nu)(\xi)}$ so $E^\alpha\in V[G*R]$.

Now $j_1(\rho_\alpha)\in j_2^*(E^\alpha)$ since

$$\forall \xi \in x_\alpha,\ f_{\kappa_1,\xi}=f_{\kappa_2,i(\xi)}\restriction \kappa_1$$
and thus $E^\alpha\in E_{\rho_\alpha}^*$.

Let us now consider $\nu_1,\nu_2\in E^\alpha$ s.t. $\pi_\alpha(\nu_1)<\pi_\alpha(\nu_2)$ and $\xi\in h_{x_\alpha}(\nu_1)$.
So we have $f_{\pi_\alpha(\nu_1),\xi}=f_{\kappa,h_{\varphi_\alpha}(\nu_1)(\xi)}\restriction \pi_\alpha(\nu_1)$.
Now
$$h_{\varphi_\alpha}(\nu_1)(\xi)=(\tilde{\sigma}_\alpha^2)^{-1}(\tilde{\sigma}_\alpha^1)^{-1}(\tau_\xi)$$
where $\tau_\xi<\pi_\alpha(\nu_1)<\pi_\alpha(\nu_2)$.
This means that
$$h_{\varphi_\alpha}(\nu_1)(\xi)\in h_{\varphi_\alpha}(\nu_2){''}h_{x_\alpha}(\nu_2)$$
and thus
$$h_{\varphi_\alpha}(\nu_2)^{-1}(h_{\varphi_\alpha}(\nu_1)(\xi))\in h_{x_\alpha}(\nu_2).$$
So we have that
$$f_{\pi_\alpha(\nu_2),h_{\varphi_\alpha}(\nu_2)^{-1}(h_{\varphi_\alpha}(\nu_1)(\xi))}=f_{\kappa,h_{\varphi_\alpha}(\nu_1)(\xi)}\restriction\pi_\alpha(\nu_2).$$
Finally we get that
$$f_{\pi_\alpha(\nu_1),\xi}=f_{\kappa,h_{\varphi_\alpha}(\nu_1)(\xi)}\restriction \pi_\alpha(\nu_1)=f_{\pi_\alpha(\nu_2),h_{\varphi_\alpha}(\nu_2)^{-1}(h_{\varphi_\alpha}(\nu_1)(\xi))}\restriction\pi_\alpha(\nu_1).$$

Let $h_{m_\alpha^*}^1$ represent $m_\alpha^*$ in $M_1$.
So in $M_2$,
$$i(m_\alpha^*)=i(j_1(h^1_{m_\alpha^*})(\rho_\alpha)
)=j_2(h^1_{m_\alpha^*})(\rho_\alpha).$$
Now $\text{dom}(m_\alpha^*)\subseteq \kappa_1 \times x_\alpha'$, and $x_\alpha'$ is determined by $\eta_\alpha^0,\ \eta_\alpha^1$ and $\eta_\alpha^2$.
Similarly, $i(x_\alpha')$ is determined by $i(\eta_\alpha^1),\ i(\eta_\alpha^2)$ and $\eta_\alpha^0$ (since $i(\eta_\alpha^0)=\eta_\alpha^0)$.
In addition,
$$\forall \xi\in x_\alpha'. i(m_\alpha^*\restriction \kappa_1\times \{\xi\})=i(m_\alpha^*)\restriction \kappa_2 \times \{i(\xi)\}.$$

We can say a bit more than that, since
$$\pi_2''\text{dom}(i(m_\alpha^*))=i(x_\alpha')=i{''}x_\alpha',$$
which means that we can compute $i(m_\alpha^*)$ from $x_\alpha'$ and $m_\alpha^*$ and vice versa.

For a condition $p$, let us denote by  $p(\beta)$ the restriction of $p$ to its $\beta$ function, i.e.
$$p(\beta):=p\restriction \kappa_i\times \{\beta\},\ i\in \{0,1,2\}.$$
First notice that since $j_1(h_{m_\alpha^*}^1)(\rho_\alpha)=m_\alpha^*$ then $j_2(h^1_{m_\alpha^*})(\rho_\alpha)=i(m_\alpha^*)$.
So we have that
$$M_2 \models \forall \xi\in x_\alpha. j_2(h_{m_\alpha^*})(\rho_\alpha,j_1(\rho_\alpha))(\xi)= j_2(h_{m_\alpha^*}^1)(\rho_\alpha)(\varphi_\alpha(\xi)).$$
So this gives that
\begin{align*}
    T:=\{(\nu_1,\nu_2)|\forall \xi\in h_{x_\alpha}(\nu_2). h_{m_\alpha^*}(\nu_1,\nu_2)(\xi)=h^1_{m_\alpha^*}(\nu_1)(h_{\varphi_\alpha}(\nu_2)(\xi))\}
    \in (E_{\rho_\alpha})^2
    \end{align*}
    and recall that
    $$h_{\varphi_\alpha}(\nu_2)(\xi)=\tilde{\sigma}_\alpha^2(\tilde{\sigma}_\alpha^1(h_{\sigma_\alpha^1}^{-1}(\nu_2)(h_{\sigma_\alpha^2}^{-1}(\nu_2)(\xi))).$$

We also have that
$$E_\alpha^\alpha:=\{\nu<\kappa|\forall \tau<h^1_{\eta_\alpha^0}(\nu).h^1_{m_\alpha^*}(\nu)(\tilde{\sigma}_\alpha^2(\tilde{\sigma}_\alpha^1(\tau)))\subseteq f_{\kappa,\tilde{\sigma}_\alpha^2(\tilde{\sigma}_\alpha^1(\tau))}\}\in E_{\rho_\alpha}.$$

Now let $T'\in E_{\rho_\alpha}$ s.t. $[T']^{*2}\subseteq T$ and let
\begin{align*}
    A'_{m^\alpha}&:=\{\nu<\kappa|f_{\kappa,\kappa^++\rho_\alpha}(\nu)\in T'\cap E_\alpha\cap E_\alpha^\alpha\},\\
    C_{m^\alpha}&:= \{\nu<\kappa|f_{\kappa,\kappa^++\rho_\alpha}(\nu)<\nu^{++}\text{ and } \nu \text{ is inaccessible}\}.
\end{align*}
Recall that $D_\alpha=\{\nu<\kappa|\pi_\alpha(f_{\kappa,\kappa^++\rho_\alpha}(\nu))=\nu\}$ and notice that $A_{m^\alpha}',C_{m^\alpha}\in U^*\cap V[G*R]$.
So we have that
$$f_{\kappa,\kappa^++\rho_\alpha}''[A'_{m^\alpha}\cap C_{m^\alpha}\cap D_\alpha]^2\subseteq A_{m^\alpha}^*$$
since for every $\nu_1<\nu_2\in A'_{m^\alpha}\cap C_{m^\alpha}\cap D_\alpha$, $f_{\kappa,\kappa^++\rho_\alpha}(\nu_i)\in T'\cap E_\alpha\cap E_\alpha^\alpha$ and since when we apply $\pi_\alpha$ we get back $\nu_1,\nu_2$, we have that
$$(f_{\kappa,\kappa^++\rho_\alpha}(\nu_1),f_{\kappa,\kappa^++\rho_\alpha}(\nu_2))\in T$$
and

$$\forall \xi\in h_{x_\alpha}(f_{\kappa,\kappa^++\rho_\alpha}(\nu_2)).f_{\kappa,h_{\varphi_\alpha}(f_{\kappa,\kappa^++\rho_\alpha}(\nu_2))(\xi)}\restriction \nu_2=f_{\nu_2,\xi}$$
from the definition of $E_\alpha$.
Combining it all together we get that
$$[A_{m^\alpha}'\cap C_{m^\alpha}\cap D_\alpha]^2\subseteq \pi_\alpha''A_{m^\alpha}^*$$
and $X_{A_{m^\alpha}}:=A_{m^\alpha}'\cap C_{m^\alpha}\cap D_\alpha$ is as required.
 \\
 $\square$

Combining the above claims, we can use Lemma \ref{the idea} to get that Gal($\kappa,W,\kappa^{++}$), which finishes the proof.

So, $W$ is a $\kappa$-complete ultrafilter  which does not satisfy the Galvin property at $\kappa^+$, but satisfies it at $\kappa^{++}$.

\subsection{ An approach via filters }

Let present now a somewhat different approach in order to show Gal$(\kappa, W, \kappa^{++})$.

Let $R=\l f_{\kappa\alpha}\mid \kappa^+\leq \alpha<\kappa^{++}\r$, i.e. the set of generic over $V[G]$ Cohen function with idexes in the interval $[\kappa^+,\kappa^{++}).$

Work in $V[G,R]$ and define $U'\subseteq \calP(\kappa)$ as follows:

$$A\in U' \text{ iff } \exists p\in G*R\quad \exists \xi \quad M_1\models(p^\frown p_\xi \Vdash_{j_1(P_{\kappa+1})} \kappa\in j_1(\lusim{A})),$$
where $\l p_\xi\mid \xi\r$ is the master condition sequence for the forcing above $\kappa$ and at $\kappa_1$ but there it is restricted to Cohens in the interval $[j_1(\kappa^+),j_1(\kappa^{++}))$, i.e. without $
Cohen(\kappa_1, \kappa_1^+)$.

\begin{lemma}
$U'$ is a normal filter over $\kappa$ in $V[G,R]$.

\end{lemma}
\pr
Clearly, $U'$ is a $\kappa-$complete filter  in $V[G,R]$.
Let us argue that it is normal. Thus let $\{A_\alpha\mid \alpha<\kappa\}\subseteq U'$ and $A=\Delta_{\alpha<\kappa} A_\alpha$.
For every $\alpha<\kappa$ pick $p^\alpha\in G*R$ and $\xi^\alpha$ witnessing that $A_\alpha\in U'$.
Let $\xi^*=\sup_{\alpha<\kappa}\xi^\alpha$.
Suppose that there is no $p\in G*R$ such that $p^\frown p_{\xi^*} \Vdash \kappa\in j_1(\lusim{A})$.
Then there are $\beta<\kappa$, $q\in G*g_\kappa$ and $\xi\geq \xi^*$ such that
$$q^\frown p_\xi\Vdash \kappa\not \in j_1(\lusim{A}_\beta).$$
But this is impossible since $q$ is compatible with $p_\beta$ and $p_\xi\geq p_{\xi_\beta}$.
\\
$\square$

The next lemma follows from the definitions.

\begin{lemma}\label{lem-hom}
$U'\subseteq U^*\cap V[G,R]$.

\end{lemma}

\begin{remark}
Note that since $U^*$ is an ultrafilter in the full generic extension, $ U^*\cap V[G,R]$ is an ultrafilter over  $V[G,R]$. However,  $U^*\cap V[G,R]\not \in   V[G,R]$, since otherwise each
Cohen function $f_{\kappa \xi} $ will be in $ V[G,R]$ as well. Namely, for every $\xi<\kappa^+$, $\tau,\rho<\kappa$,
$$f_{\kappa \xi}(\tau)=\rho \text{ iff } \{ \nu<\kappa\mid f_{\nu h_\xi(\nu)}(\tau)=\rho\} \in U^*,$$
where $h_\xi\in V$ is the canonical function which represents $\xi$.

\end{remark}

Let us define now a two dimensional version of $\tilde{U}$ of $U'$.
Let $A\in V[G*R]$ be subset of $[\kappa]^2$.
Set
$$A\in \tilde{U}   \text{ iff } \exists p\in G*R\quad \exists \xi \quad \text{ in } M_2 \quad (p^\frown p_\xi{}^\frown p^1_\xi \Vdash (\kappa,\kappa_1)\in j_2(\lusim{A})),$$
where $\l p_\xi^1\mid \xi\r$ is the master condition sequence for the forcing above $\kappa_1$ which is the image of $\l p_\xi\mid \xi\r$ under the connecting embedding $j_{12}$ from $M_1$ to $M_2$.

We have the following analog of \ref{lem-hom}:

\begin{lemma}\label{lem-hom2}
$\tilde{U}\subseteq (U^*)^2\cap V[G,R]$.

\end{lemma}

Let $A\in \tilde{U}$.
Then there are $p \in G*R$ and $\xi$ such that, in $M_2$,

$$p^\frown p_\xi{}^\frown j_{12}(p_\xi) \Vdash (\kappa,\kappa_1)\in j_2(\lusim{A}).$$

We have $j_2=j_{12}\circ j_1$.

For $\nu<\kappa$, let us denote $A_\nu=\{ \rho<\kappa\mid (\nu, \rho) \in A\}$.

Then, in $M_2$,

$$p^\frown p_\xi{}^\frown j_{12}(p_\xi) \Vdash \kappa_1\in j_{12}(j_1(\lusim{A})_\kappa).$$

Let $h_{p_\xi}$ be a function that represents $p_\xi$ in $M_1$, i.e. $j_1(h_{p_\xi})(\kappa)=p_\xi$.

Consider the set

$$X=\{\nu<\kappa\mid \text{in } M_1, p^\frown h_{p_\xi}(\nu){}^\frown p_\xi \Vdash \kappa\in j_{1}(\lusim{A}_\nu)\}.$$

Note that for every $\nu<\kappa$, if $h_{p_\xi}(\nu)\in G*R$, then $A_\nu\in U'$.

Also, the set $\{\nu<\kappa\mid h_{p_\xi}(\nu)\in G*R\}$ is in $U'$, as witnessed by $p_\xi$.

Hence,

$$Y=\{\nu<\kappa\mid A_\nu\in U'\}\in U'.$$

Set $B=Y\cap\Delta_{\nu\in Y}A_\nu$. By normality, $B\in U'$.

\begin{lemma}\label{lem5}
$[B]^2\subseteq A$.

\end{lemma}
\pr
Let $(\nu,\rho)\in [B]^2$. Then $\nu<\rho$, both in $Y$, so, $\rho\in A_\nu$. Hence, $(\nu, \rho)\in A$.
\\
$\square$

It follows now:

\begin{lemma}

$(U')^2=\tilde{U}$.

\end{lemma}
\pr
By \ref{lem5}, $(U')^2\supseteq\tilde{U}$.
The opposite inclusion follows from the definitions of $U'$ and $\tilde{U}$.
\\
$\square$

Finally recall that by \ref{lem-hom2}, $\tilde{U}\subseteq (U^*)^2\cap V[G,R]$. Hence, $$(U')^2=\tilde{U}\subseteq (U^*)^2\cap V[G,R]=W\cap V[G,R].$$
So we can run the Galvin argument with $U'$ inside $V[G,R]$ and it will give the desired conclusion for $W$.

Note that the missing parts of master condition sequences used to define $U'$ and $\tilde{U}$, i.e., $Cohen(\kappa,\kappa^+)$ and $Cohen(\kappa_1,(\kappa_1^+)^{M_1})$ are of cardinality $\kappa^+$ for the former or the image
from the ultrapower by the normal $U$ for the later, which is also of cardinality $\kappa^+$.
So, we can stabilize the needed conditions from this forcings.

\newpage
\section{ Some generalizations}
We showed above that it is consistent to have a $\kappa$-complete ultrafilter $W$ which satisfies the Galvin property at $\kappa^{++}$ but not at $\kappa^+$.
With any given cardinal $\lambda\geq \kappa$, the same method, with minor changes, will also work to get a model in which $2^\kappa=\lambda^+$ and there is a $\kappa$-complete ultrafilter $W$ which satisfies the Galvin property at $\lambda^+$ but not at $\lambda$ (assuming, say, that $\kappa$ is   $\lambda^+-$strong).
This however does not cover all options.
What if we want the first stage where the Galvin property holds to be a limit cardinal?
In this case our proof will not work anymore.
Moreover, it turns out that not everything is possible.

\begin{proposition}\label{cflesskappa}
Assume that $\kappa$ is measurable, $\lambda,\kappa<\lambda\leq 2^{\kappa},$ is a limit cardinal with cf$(\lambda)<\kappa$ and that $\lnot$Gal$(\kappa,W,\alpha)$ holds for every $\alpha<\lambda$.
Then $\lnot$Gal$(\kappa,W,\lambda)$.
\end{proposition}
\pr
    Let $\langle A^\alpha_i | i<\alpha \rangle$ witness that $\lnot$Gal$(\kappa,W,\alpha)$ $\forall\alpha<\lambda$ and let $\langle \alpha_j|j<\text{cf}(\lambda)\rangle$ be cofinal in $\lambda$.
    Assume by contradiction that Gal$(\kappa,W,\lambda)$ and look at
    $$\langle A_i^{\alpha_j}|j<\text{cf}(\lambda)\rangle.$$
    Since this set has size $\lambda$, there are $\kappa$ many sets that their intersection belongs to $W$.
    But then, since cf$(\lambda)<\kappa$, there is some $\alpha$ s.t.
     $\kappa$ many of them are from $\langle A_i^\alpha | i<\alpha\rangle$ and their intersection is in $W$, contradicting the fact that it is witnessing $\lnot$Gal$(\kappa,W,\alpha)$.
 \\
 $\square$

\begin{proposition}
    \label{cfeqkappa}
Assume that $\kappa$ is measurable, $\lambda,\kappa<\lambda\leq 2^{\kappa},$ is a limit cardinal with cf$(\lambda)=\kappa$ and that $\lnot$Gal$(\kappa,W,\alpha)$ holds for every $\alpha<\lambda$.
Then $\lnot$Gal$(\kappa,W,\lambda)$.
\end{proposition}
\pr
        Let $\langle A^\alpha_i | i<\alpha \rangle$ witness that $\lnot$Gal$(\kappa,W,\alpha)$ $\forall\alpha<\lambda$ and let $\langle \alpha_j|j<\kappa\rangle$ be cofinal in $\lambda$.
    Assume by contradiction that Gal$(\kappa,W,\lambda)$ and look at
    $$\langle B_i^{\alpha_j}|j<\kappa,i<\alpha_j \rangle \text{ where } B_i^{\alpha_j}=A_i^{\alpha_j}\setminus j.$$
    Notice that $\langle B_i^{\alpha_j}|i<\alpha_j\rangle$ is still a witness for $\lnot$Gal$(\kappa,W,\alpha_j)$.
    Since $\langle B_i^{\alpha_j}|j<\kappa,i<\alpha_j \rangle$ has size $\lambda$, there is some
    $$F\subseteq \langle B_i^{\alpha_j}|j<\kappa,i<\alpha_j \rangle, \ |F|=\kappa$$
    with $\cap F\in W$.
    If there is some $j^*$ s.t. $\kappa$ many sets in $F$ belong to $\langle B_i^{\alpha_{j^*}}|i<\alpha_{j^*}\rangle$, then we get a contradiction to the fact that $\langle B_i^{\alpha_{j^*}}|i<\alpha_{j^*}\rangle$ witness $\lnot $Gal($\kappa,W,\alpha_{j^*}$).
    Otherwise, there is some $F'\subseteq F$ and an increasing sequence $\langle j_\delta|\delta<\kappa\rangle$ s.t.
    $$F'=\{B_{i(\delta)}^{\alpha_{j_\delta}}|\delta<\kappa\} \text{ and } \cap F'\in W.$$
    But since $B_{i(\delta)}^{\alpha_{j_\delta}}\cap  j_\delta=\emptyset$, we get that $F'\cap  j_\delta=\emptyset$ for every $\delta<\kappa$, so $\cap F'=\emptyset$ which is a contradiction.
 \\
 $\square$

\subsection{$\aleph_{\kappa^{++}}$ and beyond }

We do not know whether $\aleph_{\kappa^+}$  can be a breaking point, i.e., wether it is possible to have a $\kappa-$complete ultrafilter $W$ over $\kappa$ such that

\begin{enumerate}
  \item $2^\kappa\geq \aleph_{\kappa^+}$,
  \item $\neg$Gal$(\kappa, W, \alpha)$, for every $\alpha<\aleph_{\kappa^+}$,
  \item Gal$(\kappa,W, \aleph_{\kappa^+})$.
\end{enumerate}

However, it is possible at $\aleph_{\kappa^{++}}$ or at any $\lambda$ of cofinality $\geq \kappa^{++}$.
Let us sketch the argument.

Let $\lambda=\aleph_{\kappa^{++}}$.
\\Proceed as in the construction for $\kappa^{++}$ above.  We have Cohen functions $\l f_{\kappa \alpha} \mid \alpha<\lambda \r$.
Let $$A_\alpha=\{ \nu<\kappa \mid f_{\kappa \alpha}(\nu)=1 \text{ (or odd) }\}.$$

We define a normal ultrafilter $U^*$ a normal ultrafilter over $\kappa$ as before. Then a $\kappa-$complete ultrafilter $W$ will be defined.

Let us arrange non-Galviness.

For every $\xi<\kappa^{++}$, we change values of functions $f_{\kappa_2, \gamma}(\kappa_1)$ for $\gamma$'s in the interval $(j_2(\kappa)^{+j_2(\xi)}, j_2(\kappa)^{+j_2(\xi+1)})$
if $\gamma$ has a pre-image under $k$ but not under $j_2$. This, by the usual argument, will insure $\neg Gal(\kappa, \kappa^{+\xi+1})$
using $\{A_\zeta \mid \kappa^{+\xi}<\zeta<\kappa^{+\xi+1}\}$.

Note, and this is crucial, that nothing is done in intervals of the form $(j_2(\kappa)^{+\tau}, j_2(\kappa)^{+\tau+1})$ with $\tau=\kappa$ or cofinality $\kappa$.

In particular, we cannot use  $\{A_\zeta \mid \zeta<\kappa^{+\kappa}\}$ to witness $\neg Gal(\kappa, \kappa^{+\kappa})$.
Just, for example, the sequence $\l A_{\kappa^{+\mu}} \mid \mu<\kappa \r$ may be problematic, since no changes are done inside the interval
$(j_2(\kappa)^{+\kappa}, j_2(\kappa)^{+j_2(\kappa)})$.

Turn now to Gal$(\kappa,W, \lambda)$, where $W$ is result of changes above.

Let $\{ B_\rho \mid \rho<\lambda\}\subseteq W$. Denote by $b_\rho\subseteq \lambda, |b_\rho|\leq \kappa$ the support of $B_\rho$.
Shrink to $\rho_\tau$'s with $\bigcup_{\tau<\kappa^{++}} b_{\rho_\tau}$ unbounded in $\lambda$.
Now we form a $\Delta-$system with the union still unbounded in $\lambda$, by shrink the family if necessary.
\\ The rest of the argument is similar to $\kappa^{++}$ case.

\end{document}